\documentclass[11pt]{article}
\usepackage{amsfonts}
\usepackage{mathrsfs}
\usepackage{amsfonts}
\usepackage{latexsym,amssymb,amsmath}
\usepackage{latexsym}
\usepackage{amssymb}
\usepackage{amsmath}
\usepackage{amsthm}
\usepackage{indentfirst}
\usepackage{color}
\usepackage{graphicx}
\usepackage{bm}
\usepackage{enumerate}
\usepackage{multicol}
\usepackage[square, comma, sort&compress, numbers]{natbib}
\usepackage[colorlinks]{hyperref}
\numberwithin{equation}{section} \allowdisplaybreaks

\makeatletter
\renewcommand\@cite[1]{#1\hspace{0.2em}}
\makeatother
\newtheorem{theorem}{\color{black}\indent Theorem}[section]
\newtheorem{lemma}[theorem]{\color{black}\indent Lemma}
\newtheorem{proposition}[theorem]{\color{black}\indent Proposition}
\newtheorem{definition}[theorem]{\color{black}\indent Definition}
\newtheorem{remark}[theorem]{\color{black}\indent Remark}
\newtheorem{corollary}[theorem]{\color{black}\indent Corollary}

\usepackage{amssymb,amsmath}
\DeclareMathOperator{\dive}{div}
\def\rd{{\rm d}}
\def\re{{\rm e}}

\begin{document}	
	\title{\LARGE\bf   Asymptotic behavior of the solution with positive
		temperature in nonlinear 3D thermoelasticity
	\author{$\rm{Chuang~Ma}$, $\rm{Bin~Guo}$
		\thanks{Corresponding author\newline \hspace*{6mm}{\it Email
					addresses:} machuang24@mails.jlu.edu.cn(Chuang Ma),~bguo@jlu.edu.cn~(Bin Guo)
				}
			\\
			School of Mathematics, Jilin University, Changchun 130012, PR China
	}}
	\date{ \today } \maketitle
	
{\bf Abstract:} In this paper, we study a hyperbolic-parabolic coupled system arising in nonlinear three-dimensional thermoelasticity. We establish the global well-posedness and asymptotic behavior of solutions. Our main result shows that, a thermoelastic body asymptotically converges to an equilibrium state with a uniform temperature distribution for every initial data, determined by energy conservation. The proof of the global well-posedness is divided into some steps. To begin with, we introduce an approximate problem and derive its solvability. Next, we establish a time-independent upper bound for the temperature via Moser iteration technique. Together with an estimate of gradient of entropy, we use a functional involving the Fisher information of the temperature, which enables us to handle a delicate Gronwall-type inequality, to obtain required estimates of the higher-order derivatives. Further, we prove the strict positivity of temperature by applying Moser iteration again on the negative part of the logarithm of the temperature, followed by a uniqueness argument for the weak solution. Finally, we define a dynamical system on a proper functional phase space and analyze the $\omega$-limit set for every initial data. This work provides a complete proof of the global well-posedness and the long-time behavior in the nonlinear three-dimensional thermoelasticity system.
	
\textbf{{Keywords:}} Global~existence;~Asymptotic behavior;~Nonlinear coupling term.
\thispagestyle{empty}

\section{Introduction and main results}
This article is devoted to the study of the global-in-time existence and asymptotic behavior of a thermoelastic body. The function $u\colon (0,T)\times\Omega\rightarrow\mathbb{R}^3$ denotes the displacement of the elastic material and the function $\theta\colon (0,T)\times\Omega\rightarrow\mathbb{R}$ denotes the material temperature. We consider the following nonlinear coupled problem between temperature and displacement
\begin{align} \label{equ1.1}
	\begin{cases}
		u_{tt}-\Delta u-\Delta u_{tt}=\mu\nabla \theta,~&(x,t)\in(0,T)\times\Omega,\\
		\theta_t-\Delta\theta=\mu \theta\dive(u_t),~&(x,t)\in(0,T)\times\Omega,\\
		u=0,\quad\partial_n\theta=0,~&(x,t)\in(0,T)\times\partial\Omega,\\
		u(0,x)=u_0,~u_t(0,x)=v_0,~\theta (0,x)=\theta_0> 0,&~~\quad x\in \Omega,
	\end{cases}
\end{align}
where $\mu$ is a given constant and $\Omega$ is a bounded domain in $\mathbb{R}^3$ with sufficiently smooth boundary. It is a problem which describes oscillations and the heat in a medium. This system arises a nonlinear thermoelasticity model, see \cite{Dafermos1968,Dafermos1986,RackeZheng1997,Slemrod1981} and references therein.
So far, the global well-posedness and asymptotic behavior of solutions in three-dimensional or higher-dimensional cases remain an open problem. In the present paper, we prove the existence and uniqueness of the global-in-time weak solution for the initial boundary value problem \eqref{equ1.1} and long-time behavior in three-dimensional case.

\subsection{Physical derivation}
In this section, we briefly describe the derivation of the physical foundations of the model \eqref{equ1.1}. For a more comprehensive explanation, we refer to \cite{Dafermos1968,Dafermos1986,Slemrod1981,JR2000,CTCVPDE2024} and references therein. We consider a system of partial differential equations modeling the balance of linear momentum and energy
\begin{equation}\label{equME}
	\begin{cases}
		\rho \mathbf{u}_{tt}=\nabla\cdot\sigma,\\
		e_t+\nabla\cdot\mathbf{q}=\sigma\cdot\nabla\mathbf{u}_t,
	\end{cases}
\end{equation}
where $\rho$ is the material density, $\sigma$ is the Piola-Kirchhoff stress tensor, $e$ is the internal energy, $\mathbf{q}$ is the heat flux and the state variables are the displacement vector $\mathbf{u}$ and the temperature $\theta$.
Additionally, $\sigma\cdot\nabla\mathbf{u}_t$ represents the work performed by the elastic body.
Let us introduce the Helmholtz free energy function
$$\psi:=e-\theta S,$$
where $S$ is the entropy. Next, we assume the response functions $\sigma,~e$, and $\psi$ are independent of the temperature gradient $\nabla\theta$,
namely,
$$\sigma=\sigma(\theta,\nabla\mathbf{u}),~~~e=e(\theta,\nabla\mathbf{u}),~~~\psi=\psi(\theta,\nabla\mathbf{u}).$$
As a consequence of the Clausius-Duhem inequality \cite{JR2000}, $\sigma$ and $S$ are expressed as follows
\begin{align}\label{C-Duhem}
	\sigma(\theta,\nabla\mathbf{u})=\psi_\mathbf{F}(\theta,\nabla\mathbf{u}),~~~S(\theta,\nabla\mathbf{u})=-\psi_\theta(\theta,\nabla\mathbf{u}),
\end{align}
where the subscripts denote the partial derivatives (by $\psi_\mathbf{F}$ we denote the Frech\'et derivative with respect to the second variable).
Using \eqref{C-Duhem}, we can write these relationships in the following format of Gibbs' relations
\begin{align}\label{Gibbs}
	\begin{split}
		\psi_\theta=e_\theta-S-\theta S_\theta~~~&\Longrightarrow~~~\theta S(\theta,\nabla\mathbf{u})_\theta=e(\theta,\nabla\mathbf{u})_\theta,\\
		\psi_\mathbf{F}=e_\mathbf{F}-\theta S_\mathbf{F}~~~&\Longrightarrow~~~\theta S(\theta,\nabla\mathbf{u})_\mathbf{F}=e(\theta,\nabla\mathbf{u})_\mathbf{F}-\sigma(\theta,\nabla\mathbf{u}).
	\end{split}
\end{align}
Dividing the second equation in \eqref{equME} by $\theta$ (provided that $\theta>0$), we obtain the entropy equation
\begin{equation}\label{equSt}
	S_t+\nabla\cdot\left(\frac{\mathbf{q}}{\theta}\right)=-\frac{\mathbf{q}\cdot\nabla\theta}{\theta^2},
\end{equation}
where $\nabla\cdot\left(\frac{\mathbf{q}}{\theta}\right)$ is the entropy flux and $-\frac{\mathbf{q}\cdot\nabla\theta}{\theta^2}$ is the entropy production.
Here we have used the following equations
$$\frac{\nabla\cdot\mathbf{q}}{\theta}=\nabla\cdot\left(\frac{\mathbf{q}}{\theta}\right)+\frac{\mathbf{q}\cdot\nabla\theta}{\theta^2},$$
and
\begin{align}\label{et}
	e(\theta,\nabla\mathbf{u})_t=(\psi+\theta S)_t
	=\psi_\mathbf F\nabla\mathbf{u}_t+\theta S_t
	=\sigma\cdot\nabla\mathbf u_t+\theta S_t.
\end{align}
Based on the second law of thermodynamics \cite{JR2000}, $\mathbf{q}(\theta,\nabla\theta)$ have to obey the heat conduction inequality
$$\mathbf q\cdot\nabla\theta\leqslant0,$$
such that the entropy production $-\frac{\mathbf{q}\cdot\nabla\theta}{\theta^2}\geqslant0$.
By using \eqref{Gibbs} and \eqref{et}, the equations \eqref{equME} can be transformed into
\begin{equation}\label{equndecomposed}
	\begin{cases}
		\rho\mathbf u_{tt}=\nabla\cdot\left(e_{\mathbf{F}}-\theta S_\mathbf{F}\right),\\
		e_\theta\theta_t+\nabla\cdot\mathbf{q}=-\theta S_\mathbf{F}\nabla\mathbf u_t.
	\end{cases}
\end{equation}
We shall assume that energy can be decomposed as
\begin{align*}
	e(\theta,\nabla\mathbf u)=e_1(\theta)+e_2(\nabla\mathbf u),
\end{align*}
and the entropy can also be decomposed as
\begin{align*}
	S(\theta,\nabla\mathbf u)=S_1(\theta)+S_2(\nabla\mathbf u).
\end{align*}
Inserting the above expressions into \eqref{equndecomposed}, the original system \eqref{equME} now simplifies to
\begin{equation*}
	\begin{cases}
		\rho\mathbf u_{tt}=\nabla\cdot\left(e_2(\nabla\mathbf u)_{\mathbf{F}}-\theta S_2(\nabla\mathbf u)_\mathbf{F}\right),\\
		e_1(\theta)_\theta\theta_t+\nabla\cdot\mathbf{q}=-\theta S_2(\nabla\mathbf u)_t.
	\end{cases}
\end{equation*}
Notice that the work transferred to heat is given by $-\theta S_2(\nabla\mathbf u)_t$. It is very restrictive on $\theta$ and completely determines the coupling.
The analysis is completed by requiring thermodynamical stability conditions on the energy and entropy
\begin{align}\label{ife1}
	\lim\limits_{\theta\rightarrow0^+}S_1(\theta)=-\infty,~~e_1(\theta)>0,~~e_1(\theta)_\theta>0,
\end{align}
for all $\theta>0$ and there exist constant $C_1>0$ and $C_2$ such that
\begin{equation}\label{ife2}
	e_2(\nabla \mathbf u)-C_1S_2(\nabla\mathbf u)\geqslant C_2.
\end{equation}
Further, we assume the following boundary conditions
\begin{align*}
	(\sigma\mathbf{u_t})\cdot\mathbf n=0,
	~~~\mathbf q\cdot\mathbf n=0,~~~\text{on}~(0,T)\times\partial\Omega,
\end{align*}
and introduce the Helmholtz function
\begin{align*}
	H_{C_1}(\theta,\nabla\mathbf u):=e(\theta,\nabla\mathbf u)-C_1S(\theta,\nabla\mathbf u).
\end{align*}
By multiplying \eqref{equSt} with $C_1$ and integrating over $(0,t)\times\Omega$, we get
\begin{align}\label{HC1}
	\int_0^t\int_\Omega \bigg(H_{C_1}(\theta,\nabla\mathbf u)-e(\theta,\nabla\mathbf u)\bigg)_s\rd x\rd s-C_1\int_0^t\int_\Omega\frac{\mathbf{q}\cdot\nabla\theta}{\theta^2}\rd x\rd s=0.
\end{align}
Multiplying $\eqref{equME}_1$ with $\mathbf u_t$ and summing up with $\eqref{equME}_2$, we obtain the total energy balance
$$\frac 12\int_\Omega|\mathbf u_t|^2\rd x+\int_\Omega e(\theta,\nabla\mathbf u)\rd x=\frac 12\int_\Omega|\mathbf u_t(0)|^2\rd x+\int_\Omega e\big(\theta(0),\nabla\mathbf u(0)\big)\rd x.$$
Inserting the above equality into \eqref{HC1}, we have the total dissipation balance
\begin{align}\label{balance1}
	&~~~\int_\Omega H_{C_1}(t)\rd x+\frac 12\int_\Omega|\mathbf u_t(t)|^2\rd x-C_1\int_0^T\int_\Omega\frac{\mathbf{q}\cdot\nabla\theta}{\theta^2}\rd x\rd t\notag\\
	&=\int_\Omega H_{C_1}(0)\rd x+\frac 12\int_\Omega|\mathbf u_t(0)|^2\rd x.
\end{align}
By using \eqref{Gibbs}, we notice that the condition \eqref{ife1} implies $S_1(\theta)_\theta=\frac{e_1(\theta)_\theta}{\theta}>0$. Further, we have
\begin{align*}
	\left[H_{C_1}(\theta,\nabla\mathbf u)\right]_\theta
	&=\left[e_1(\theta)-C_1S_1(\theta)\right]_\theta\\
	&=e_1(\theta)_\theta-\theta S_1(\theta)_\theta+(\theta-C_1)S_1(\theta)_\theta\\
	&=(\theta-C_1)S_1(\theta)_\theta,
\end{align*}
which means that the function $e_1(\theta)-C_1S_1(\theta)$ exists its global minimum on $\mathbb{R}^+$ at $\theta=C_1$, defined by $H_{C_1}(\theta)_{\min}$. This together with \eqref{ife2} implies the coercivity of the Helmholtz function
\begin{align*}
	H_{C_1}(\theta,\nabla\mathbf u)
	&=e(\theta,\nabla\mathbf u)-C_1S(\theta,\nabla\mathbf u)\\
	&=\left(e_1(\theta)-C_1S_1(\theta)\right)+\left(e_2(\nabla\mathbf u)-C_1S_2(\nabla\mathbf  u)\right)\\
	&\geqslant H_{C_1}(\theta)_{\min}+C2.
\end{align*}
Thus, from \eqref{balance1}, we know that the boundedness of the entropy production.
Under the assumption of the small deformation, we incorporate the rotational inertia effects through a second-order strain energy term
 $$e_2(\nabla\mathbf u)=\frac 12|\nabla\mathbf u|^2+\frac 12\partial^2_t\left(|\nabla\mathbf u|^2\right),$$
which accounts for both elastic potential energy and rotational inertia.
For the stress due to heat expansion, we use the standard linearized constitutive relation
$\theta S_2(\mathbf u)_\mathbf{F}=\zeta\theta I$, where the constant $\zeta>0$.
Observing the relationship between the deformation gradient and displacement $(\nabla\cdot\mathbf u)_\mathbf F=I$, it immediately leads to
$$S_2(\nabla\mathbf{u})=\zeta\nabla\cdot\mathbf u.$$
Here $\nabla\cdot\mathbf u$ represents the linearization of Jacobian.
Finally, we assume that the body is fixed at the boundary and the internal energy flux through the boundary is zero, namely
$$\mathbf{u}=0,~~~\partial_n\theta=0,~~~\text{on}~(0,T)\times\partial\Omega,$$
which means that the boundary of the body is rigidly clamped and the boundary $\partial\Omega$ is kept at a constant temperature $\theta_0$ in thermodynamical terms.
Supplementing with initial data, we arrive at the following model of the balance of momentum and internal energy
\begin{align*}
\begin{cases}
	\rho\mathbf{u}_{tt}=\nabla\cdot\left(\nabla\mathbf u-\nabla\mathbf{u}_{tt}-\zeta\theta I\right),~~~&\text{in}~(0,T)\times\Omega,\\
	C_V(\theta)\theta_t+\nabla\cdot\mathbf q(\theta,\nabla\theta)=-\zeta\theta\nabla\cdot\mathbf u_t,~&\text{in}~(0,T)\times\Omega,\\
	\mathbf{u}=0,~~~\partial_n\theta=0,~~~&\text{on}~(0,T)\times\partial\Omega.\\
	u(0,x)=u_0,~u_t(0,x)=v_0,~\theta (0,x)=\theta_0> 0,&\text{in}~~\Omega,
\end{cases}
\end{align*}
where $C_V(\theta):=e_1(\theta)_\theta$ is the heat capacity and the internal energy flux $\mathbf{q}(\theta,\nabla\theta)$ satisfies the Fourier law $\mathbf{q}(\theta,\nabla\theta)=-\kappa(\theta)\nabla\theta$, where $\kappa(\theta)>0$ is the heat conductivity.
Based on the Dulong-Petit law \cite{Simon2011}, we know that the heat capacity is almost constant at high temperatures for many homogeneous materials.
On the other hand, heat conductivity tends to show little change at high temperatures in many relevant materials \cite{Ross1991}.
Then the following hypothesis is reasonable
$$C_V(\theta)=1,~~~\kappa(\theta)=1.$$
Assuming further the density of elastic body $\rho\equiv1$, we obtain exactly the model \eqref{equ1.1}.
The above model is in full accordance with the laws of thermodynamics.
\subsection{A related literature overview}
In this section, we compare our findings with known literature.
To the best of our knowledge,  the existence and asymptotic behavior as time tends to infinity for one-dimensional thermoelastic systems initially studied in the 1960s. Dafermos \cite{Dafermos1968} established the first rigorous results on long-time behavior, however the initial work was related to the linear simplification.
Nonlinear elasticity theory demands rigorous computational and analytical method beyond classical linear approaches. In 1981, Slemrod \cite{Slemrod1981} considered the following more general one-dimensional nonlinear thermoelastic problem
\begin{equation}\label{equ-Slemrod}
	\begin{cases}
	u_{tt}-a(u_x,\theta)u_{xx}=-b(u_x,\theta)\theta_x,~~~&\text{on}~(0,T)\times(0,1),\\
	c(u_x,\theta)\theta_t-d(u_x,\theta)\theta_{xx}=-b(u_x,\theta)u_{xt},&\text{on}~(0,T)\times(0,1),
	\end{cases}
\end{equation}
with mixed boundary conditions (Neumann for $u$, Dirichlet for $\theta$), where $a,~b,~c,~d$ are differentiable functions such that $a,~|b|,~c,~d>0$. The author proved local existence and uniqueness of a solution with specific assumptions based on application of the contraction mapping theorem to solutions of a related linear problem. Furthermore, assuming the initial data are sufficiently small, the author also proved global existence and uniqueness of a solution based on an a priori estimate.
The Dirichlet initial boundary value problem of \eqref{equ-Slemrod} has been studied by Racke \cite{Racke1988} who established the local existence of classical smooth solutions for both bounded and unbounded domains assuming smooth data.
Subsequently, Jiang \cite{Jiang1990} established the global existence of solutions for small smooth initial data and analyzed the decay properties of classical solutions in one-dimensional nonlinear thermoelasticity.
Cie{\'s}lak et al. \cite{CT2022} extended these results. They proved that the potential non-negativity of the temperature and existence and uniqueness of the local-in-time strong solution. They also established existence of the global-in-time weak measure valued solutions by using the vanishing viscosity method \cite{ChenGQ1995}.

Regarding singularities, Dafermos and Hsiao \cite{Dafermos1986} investigated the development of singularities for a special one-dimensional nonlinear model. They rigorously demonstrated that a smooth solution for large data will blow up in finite time.
A similar problem was also considered by Hrusa and Messaoudi \cite{MH1990}. They established the existence of smooth initial data leading to finite-time singularity formation for solutions.

For higher-dimensional thermoelastic systems, Racke \cite{Racke2009} studied the classical hyperbolic-parabolic model with Fourier's law of heat conduction, establishing global well-posedness results for smooth or weak solutions and blow-up results on linear systems and nonlinear systems. He also stated that the asymptotic behavior of solutions as time tends to infinity for the linearized system. It has been shown that solutions starting close to steady states exhibit asymptotic convergence to these equilibria in \cite{Racke2009}.
Let us also mention \cite{CTCVPDE2024}, Cie{\'s}lak  et al. first proved global existence of weak solutions with defect measure which satisfies the weak-strong uniqueness for large initial data through a key reformulation using thermodynamic justified variables to obtain an equivalent system, whereby the internal energy balance is replaced with entropy balance in nonlinear three-dimensional thermoelasticity inspired by the theory of Feireisl and Novotny from \cite{FENA2009,FENA2012}.

Recently, Bies and Cie{\'s}lak \cite{CTglobal1D,BPCT2024} utilized a Fisher information-based functional \cite{CT3Dchemorepulsion}, together with the Bernis-type inequality given by Winkler \cite{Winkler2012}, to prove the existence, uniqueness and asymptotic behavior of a regular global-in-time solution for system \eqref{equ-Slemrod} with positive temperature in nonlinear one-dimensional thermoelasticity under appropriate initial conditions, where the thermodynamical parameter functions $a,~c,~d$ are identically $1$ and $b\in\mathbb{R}$ arbitrary.

Over the past four decades, there have been attempts to investigate more nonlinear models, under some particular choices of nonlinear term $p(u_{tx},\theta)$ replacing the term $-\mu\theta\dive(u_t)$ in the heat equation of \eqref{equ1.1}. However, most of these studies are limited to special one-dimensional cases or depend on particular forms of nonlinearities in the developed theory. Furthermore, none of the works mentioned above dealt with the asymptotic behavior of solutions as time tends to infinity in higher-dimensional nonlinear thermoelastic systems, which is, besides the proof of global-in-time existence, uniqueness and strict positivity of the solution, one of the main novelties of our paper.

In this paper, we study the global-in-time existence and asymptotic behavior of solutions for general initial data for nonlinear three-dimensional thermoelastic system \eqref{equ1.1}. We are interested in positive solutions, namely, $\theta(t,x)>0$ for almost all $(t,x)\in[0,\infty)\times\Omega$. Indeed, our solutions are proven to be positive. We also demonstrate that weak solutions (defined in the next section) converge with time to the particular steady state whatever the data is. Our main results are the following theorems:

\begin{theorem}\label{ThexistandPos}
	For any constant $\mu\in\mathbb{R}$, there exists a unique global-in-time solution of \eqref{equ1.1} in the sense of Definition \ref{solution} with a positive temperature.
\end{theorem}

\begin{theorem}\label{Th-Convergences}
Let us assume that $u$ and $\theta$ be a weak solution of \eqref{equ1.1} with initial data $u_0,~v_0$ and $\theta_0>0$, respectively.
Then, the following convergences
\begin{align*}
	u(t,\cdot)&\rightarrow0~~~~\text{in}~H_0^1(\Omega),~\text{when}~t\rightarrow\infty,\\
	u_t(t,\cdot)&\rightarrow0~~~~\text{in}~H_0^1(\Omega),~\text{when}~t\rightarrow\infty,\\
	\theta(t,\cdot)&\rightarrow\theta_\infty~~\text{in}~L^2(\Omega),~\text{when}~t\rightarrow\infty,
\end{align*}
hold, where $\theta_\infty=\left(\frac 12\int_\Omega v_0^2\rd x+\frac 12\int_\Omega |\nabla u_{0}|^2\rd x+\frac{1}{2}\int_\Omega|\nabla v_0|^2\rd x+\int_\Omega\theta_0 \rd x\right)\cdot|\Omega|^{-1}$.
\end{theorem}
From a physical perspective, the above theorem shows that the dissipative nature of heat transport induces decay of mechanical oscillations in thermoelastic systems. Let us point out the above results are also applicable in lower dimensions, since the Sobolev embedding theorem has better properties in lower dimensional $1$ and $2$.
Next, we outline the method of proof, highlighting how physics-driven estimates govern the analysis.

\textbf{Mathematical challenge.} The main difficulties for proving the above theorems are due to the fact that the energy equation \eqref{equ1.1} is highly nonlinear by the presence of the nonlinear terms $\mu\theta\dive(u_t)$, where the usual compactness based arguments do not apply. A key observation is the reformulation of the approximate problem as an entropy equation provided that the temperature is positive, which yields a rigorous quantitative version of the second law of thermodynamics. Since we work with weak solutions, we combine the energy method with the Moser-type iteration technique \cite{Moser} used by Alikakos \cite{Alikakos1979} to obtain a delicate estimate in $L^\infty$, instead of defect measure estimate as in \cite{CTCVPDE2024}. It makes up for the lack of the Sobolev embedding $H^1\hookrightarrow L^\infty$ in three-dimensional case.
Together with a basic lemma in analysis given in \cite{Winkler2012,CT3Dchemorepulsion}, we get the required compactness estimates to prove that the global-in-time existence and uniqueness of solutions in the sense of Definition \ref{solution}. This also makes the improvement with respect to \cite{CTCVPDE2024} possible.
To study the asymptotic behavior, we define the dynamical system on a proper functional phase space and observe its $\omega$-limit set by using the quantitative version of the second principle of thermodynamics. Furthermore, both the entropy dissipation term and the lower bound of entropy play crucial roles when identifying the structure of the $\omega$-limit set.

\textbf{Structure of the paper.}~In Section $2$, we state the definition of weak solutions to \eqref{equ1.1} and some preliminary information, including some useful lemmas to obtain the required compactness estimates.
In Section $3$, we estabilish the global-in-time existence of \eqref{equ1.1} via the approximation with the half-Galerkin method and uniqueness with positive temperature to complete the proof of Theorem \ref{ThexistandPos}. The proof be divided into some steps.
To begin with, we introduce an approximate problem \eqref{equAppro} and prove its solvability by using Schaefer's fixed point theorem \cite{Evans1}. In order to pass to the limit of \eqref{equAppro}, we still need that the approximation sequence has sufficiently strong compactness properties. The key observation is that we can rewrite approximate problem \eqref{equAppro} in the form of entropy equation \eqref{equ-entropy} provided that $\theta^n>0$, which leads to a mathematically rigorous quantification of the second principle of thermodynamics.
Next, we get a time-independent upper bound of temperature via the Moser iteration \cite{Moser} by the method of Alikakos \cite{Alikakos1979}.
With the help of an estimate of the gradient of entropy, we employ a functional from \cite{CT3Dchemorepulsion} that involves the Fisher information of the temperature. This allows us to control a delicate Gronwall-type inequality, ultimately yielding the required time-independent higher-order estimates.
Subsequently, we prove the strict positivity of temperature by applying Moser iteration again on the negative logarithm of the temperature.
Further, the uniqueness of the weak solution to \eqref{equ1.1} was obtained.
The final section is devoted to the asymptotic analysis of solutions as time tends to infinity, which completes the proof of Theorem \ref{Th-Convergences}.
The proof proceeds by defining the dynamical system on a well-defined functional phase space. Then we can identify the corresponding $\omega$-limit set for every initial data. The key technical ingredients include the quantitative second law of thermodynamics, the entropy dissipation estimate, and lower bound of entropy to identify the asymptotic state of a heat-conducting elastic body.

\section{Preliminaries}
Throughout this paper, $C$ is a positive constant and may be different at each occurrence.
Before proving our main results, we quote some lemmas which will be used in our arguments.
\begin{lemma}[Aubin-Lions]\cite{Simon1986}\label{LemAL}
Let $X,~Y$ and $Z$ be the Banach spaces, $X\hookrightarrow\hookrightarrow Y\hookrightarrow Z$,
\begin{align*}
	W_1&=\left\{u\in L^p(0,T; X);~u_t\in L^1(0,T; Z)\right\}~\text{with}~1\leqslant p<\infty,\\
	W_2&=\left\{u\in L^\infty(0,T; X);~u_t\in L^r(0,T; Z)\right\}~\text{with}~r>1.
\end{align*}
Then
$$W_1\hookrightarrow\hookrightarrow L^p(0,T; Y)~~~\text{and}~~~W_2\hookrightarrow\hookrightarrow C([0,T]; Y).$$
The sign $\hookrightarrow\hookrightarrow$ and $\hookrightarrow$ denote compact embedding and continuous embedding, respectively.
\end{lemma}

The following three-dimensional variant of the functional inequality derived in \cite{CT3Dchemorepulsion} provides a key tool to control the Fisher information dissipation of the temperature in our framework.
\begin{lemma}\cite{Winkler2012,CT3Dchemorepulsion} \label{LemDsqrt}
Let $\Omega\subset\mathbb{R}^3$ be a smooth bounded domain. For every positive function $\varphi\in C^2(\overline{\Omega})$ with the boundary condition $\frac{\partial \varphi}{\partial n}\big|_{\partial\Omega}=0$, we have
\begin{align*}
	\int_\Omega|D^2\sqrt{\varphi}|^2\rd x
	\leqslant \frac{11+4\sqrt{3}}{8}\int_\Omega \varphi|D^2\log \varphi|^2\rd x.
\end{align*}
\end{lemma}

In what follows, we give the definition of weak solutions of \eqref{equ1.1}.
\begin{definition}\label{solution}
We say that $u,~\theta$ is a solution to \eqref{equ1.1} if the following hold:
\begin{enumerate}[$\bullet$]
	\item The initial data is of regularity
	$$u_0\in H^2(\Omega)\cap H_0^1(\Omega),~~~v_0\in H^2(\Omega)\cap H_0^1(\Omega),~~~\theta_0\in H^1(\Omega).$$
	We also require that there exists $\tilde\Theta>0$ such that $\theta_0\geqslant\tilde\Theta$ for all $x\in\Omega$.
	\item Solutions $\theta$ and $u$ satisfy
	\begin{align*}
		&\theta\in L_{\rm loc}^\infty\left(0,\infty; H^1(\Omega)\right)\cap H_{\rm loc}^1\left(0,\infty; L^2(\Omega)\right),\\
		&u\in L_{\rm loc}^\infty\left(0,\infty; H^2(\Omega)\right)\cap W_{\rm loc}^{1,\infty}\left(0,\infty; H_0^1(\Omega)\cap H^2(\Omega)\right)\cap W_{\rm loc}^{2,\infty}\left(0,\infty; H_0^1(\Omega)\right).
	\end{align*}
	\item The momentum equation
	\begin{align*}
		\int_\Omega u_{tt}\varphi \rd x+\int_\Omega \nabla u\nabla\varphi \rd x+\int_\Omega\nabla u_{tt}\nabla\varphi\rd x
		=-\mu \int_\Omega\theta\nabla\cdot\varphi\rd x
	\end{align*}
	is satisfied for all $\varphi\in C^\infty(\Omega)$ and for almost all $t\in(0,\infty)$.
	\item The entropy equation
	$$\int_\Omega\theta_t\psi\rd x+\int_\Omega\nabla\theta\cdot\nabla\psi\rd x=\mu\int_\Omega\theta\dive(u_t)\psi\rd x$$
	holds for all $\psi\in C^\infty(\Omega)$ and for almost all $t\in(0,\infty)$.
	\item Initial conditions are attained in the following sense:
	$$\theta\in C\left([0,\infty); L^2(\Omega)\right),$$
	$$u\in C\left([0,\infty); H_0^1(\Omega)\right)~~~\text{and}~~~u_t\in C\left([0,\infty); H_0^1(\Omega)\right).$$
\end{enumerate}
\end{definition}

Last but not least, the following energy conservation identity of \eqref{equ1.1} is obtained.
\begin{proposition}Regular solutions $u,~\theta$ of \eqref{equ1.1} satisfy
	\begin{align} \label{energy}
		&~~~\frac 12\int_\Omega u_t^2\rd x+\frac 12\int_\Omega |\nabla u|^2\rd x+\frac{1}{2}\int_\Omega|\nabla u_t|^2\rd x+\int_\Omega\theta \rd x\notag\\
		&=\frac 12\int_\Omega v_0^2\rd x+\frac 12\int_\Omega |\nabla u_{0}|^2\rd x+\frac{1}{2}\int_\Omega|\nabla v_0|^2\rd x+\int_\Omega\theta_0 \rd x:=E_0.
	\end{align}
\end{proposition}
\begin{proof}
Multiplying the first equation of \eqref{equ1.1} by $u_t$ and integrating over $\Omega$, it yields
	\begin{align}\label{pro1}
		\frac 12\frac{\rd}{\rd t}\left(\int_\Omega u_t^2\rd x+\int_\Omega |\nabla u|^2\rd x+\int_\Omega|\nabla u_t|^2\rd x\right)=\mu \int_\Omega\nabla\theta \cdot u_t \rd x.
	\end{align}
Integrating the second equation of \eqref{equ1.1} over $\Omega$ directly, we get
   \begin{align} \label{pro2}
	\frac{{\rd}}{\rd t}\int_\Omega\theta \rd x=-\mu \int_\Omega\nabla\theta u_t\rd x.
   \end{align}
Now we add \eqref{pro2} to \eqref{pro1} to obtain
\begin{align*}
	\frac{{\rd}}{\rd t}\left(\frac 12\int_\Omega u_t^2\rd x+\frac 12\int_\Omega |\nabla u|^2\rd x+\frac 12\int_\Omega|\nabla u_t|^2\rd x+\int_\Omega\theta \rd x\right)=0.
\end{align*}
Therefore the total energy balance \eqref{energy} is obtained.
\end{proof}

\section{Proof of Theorem \ref{ThexistandPos}}
This section is devoted to the proof of Theorem \ref{ThexistandPos}. Next, we will be devote to the rigorous construction of a global-in-time regular solution to \eqref{equ1.1} such that the temperature $\theta$ is strict positive. More precisely, in the wave equation part of \eqref{equ1.1} we will use the Galerkin approximation method. The heat equation part of \eqref{equ1.1} will be solved in a straight manner.

\subsection{Approximate problem}
As a rule, there are various ways to construct weak solutions. It seems to us that the half-Galerkin method we have selected achieves the desired goal most efficiently.
Let us introduce a smooth orthonormal basis  $\{\varphi_i\}$ of $H_0^1(\Omega)$ and denote $V_n=\operatorname{span}\{(\varphi_i,0,0), (0,\varphi_i,0), (0,0,\varphi_i)\}_{1\leqslant i\leqslant n}$.
\begin{definition}\label{definitionApproximate}
We say that $u_n\in W_{\rm loc}^{2,\infty}\left(0,\infty; V_n\right)$ and $\theta_n\in L_{\rm loc}^{2}\left(0,\infty; H^1(\Omega)\right)$ is a solution of an approximate problem if for all $\varphi\in V_n$ and $\psi \in H^1(\Omega)$ the following equations are satisfied in $\mathcal{D}'(0,\infty)$
\begin{align} \label{equAppro}
 	\begin{split}
 		&\int_\Omega u^n_{tt}\varphi \rd x+\int_\Omega \nabla u^n\nabla\varphi \rd x+\int_\Omega\nabla u_{tt}^n\nabla\varphi\rd x=-\mu \int_\Omega\theta^n\nabla\cdot\varphi\rd x,\\
 		&\int_\Omega\theta^n_t\psi \rd x+\int_\Omega \nabla\theta^n\nabla\psi\rd x=\mu \int_\Omega\theta^n(\nabla\cdot u_t^n)\psi\rd x.
 	\end{split}
\end{align}
Moreover, the following initial equalities
$$\theta^n(0)=\theta_0^n,~~~u^n(0)=P_{V_n}u_0,~~~u_t^n(0)=P_{V_n}v_0$$
hold, where $P_{V_n}$ is an orthogonal projection of the corresponding spaces onto $V_n$.
Moreover, for the initial data, we choose $0<c_n\leqslant\theta_0^n\in H^1(\Omega)$ as a regularization of $\theta_0$ such that $\partial_n\theta_0^n=0$ on $\partial\Omega$,~$\theta_0^n\rightarrow\theta_0$ in $L^1(\Omega)$ as $n\rightarrow\infty$ and $\int_\Omega\theta_0^n\rd x\leqslant\int_\Omega\theta_0\rd x$ for all $n\in\mathbb{N}$.
\end{definition}

\begin{lemma}\label{Lemappenergy}
	Let $u_n\in W_{\rm loc}^{2,\infty}\left(0,\infty; V_n\right)$ and $\theta_n\in L_{\rm loc}^{2}\left(0,\infty; H^1(\Omega)\right)$ is a solution to the approximate problem given by Definition \ref{definitionApproximate}.
	Then the following energy equality is satisfied:
	\begin{align*}
		&~~~\frac 12\int_\Omega |u_t^n|^2\rd x+\frac 12\int_\Omega |\nabla u^n|^2\rd x+\frac{1}{2}\int_\Omega|\nabla u_t^n|^2\rd x+\int_\Omega\theta^n \rd x=E_0^n,
	\end{align*}
where $E_0^n=\frac 12\int_\Omega |P_{V_n}v_0|^2\rd x+\frac 12\int_\Omega |\nabla P_{V_n}u_{0}|^2\rd x+\frac{1}{2}\int_\Omega|\nabla P_{V_n}v_0|^2\rd x+\int_\Omega\theta_0^n\rd x$ represents the energy of the initial data as projected onto the finite-dimensional subspace $V_n$.
\end{lemma}
\begin{proof}
	Let us take $\varphi=u^n_t$ in $\eqref{equAppro}_1$ and $\psi=1$ in $\eqref{equAppro}_2$ and add the resulting equations to obtain
	\begin{align*}
			\frac{\rd}{\rd t}\left(\frac 12\int_\Omega |u_t^n|^2\rd x+\frac 12\int_\Omega |\nabla u^n|^2\rd x+\frac{1}{2}\int_\Omega|\nabla u_t^n|^2\rd x+\int_\Omega\theta^n \rd x\right)=0.
	\end{align*}
Further, we integrate over $(0,t)$. Thus, the claim is shown.
\end{proof}
\begin{remark}
	Notice that $\lim\limits_{n\rightarrow\infty}E_0^n=E_0$ where $E_0$ defined in \eqref{energy} is independent of time variable, and $E_0^n\leqslant E_0$ for all $n\in\mathbb{N}$. This will be very useful in the next section when we pass to the limit for approximation problem in sense of Definition \ref{definitionApproximate}.
\end{remark}

The strict positivity of temperature $\theta^n(t,x)$ enters into Proposition \ref{thetangeq0} through the following Lemma from \cite[Lemma 4.2]{CTCVPDE2024}; we do not reproduce the proof here.
\begin{lemma} \label{theta+} 
	Let $f\in L^1\left(0,T; L^\infty(\Omega)\right)$ and $\theta\in L^2\left(0,T; H^2(\Omega)\right)\cap H^1\left(0,T; L^2(\Omega)\right)$ be a solution to the following parabolic problem
	\begin{equation*}
		\begin{cases}
			\theta_t-\Delta \theta=f\cdot\theta,~~~&(t,x)\in(0,T)\times\Omega,\\
			\partial_n\theta=0,&(t,x)\in(0,T)\times\partial\Omega,\\
			\theta(0,\cdot)=\theta_0^n>0.
		\end{cases}
	\end{equation*}
Then we have $\theta(t,x)>0$ for all $x\in\Omega$ and any $t\geqslant0$.
\end{lemma}

\begin{proposition}\label{thetangeq0}
	For every $n\in\mathbb{N}$, there exists a solution $(u^n,\theta^n)$ of the approximate problem in the sense of Definition \ref{definitionApproximate}. Moreover one has $\theta^n\in L_{\rm loc}^2\left(0,\infty; H^2(\Omega)\right)\cap H_{\rm loc}^1\left(0,\infty; L^2(\Omega)\right)$ and $\theta^n(t,x)>0$ for all $x\in\Omega$ and any $t\geqslant0$.
\end{proposition}
\begin{proof}
We prove the existence to the approximate problem in the sense of Definition \ref{definitionApproximate} by using Schaefer's fixed point theorem \cite{Evans1}. Next, we divide the proof into three steps.

\textbf{Step 1.~Denote the fixed point mapping $\mathcal{A}$}. Let us choose an arbitrary $n\in\mathbb{N}$ (but fixed) and introduce the fixed point set
$$X_n=C^1\left([0,T]; V_n\right).$$
We equip $X_n$ with the following norm
$$\|u_n\|_{X_n}=\max\limits_{0\leqslant t\leqslant T}\left\{\|u(t)\|_{L^2(\Omega)},~\|u_t(t)\|_{L^2(\Omega)}\right\}.$$
Next, we define the fixed point mapping $\mathcal{A}: X_n\rightarrow X_n$ in the following way.
Let $\mathbf u\in X_n$ and define $\theta(\mathbf u)\in X_n$ as a solution of the following initial boundary value problem
\begin{equation*}
	\begin{cases}
		\theta(\mathbf u)_t-\Delta\theta(\mathbf u)=\mu \theta(\mathbf u)\dive(\mathbf u_t),~~~&\text{in}~\Omega\times(0,T),\\
		\partial_n\theta(\mathbf u)=0,~&\text{on}~\partial\Omega\times(0,T),\\
		\theta (0,\cdot )=\theta_0> 0.
	\end{cases}
\end{equation*}
More precisely, $\theta(\mathbf u)$ satisfies the following equation
\begin{equation}\label{theta-mathbf-u}
\int_\Omega\theta(\mathbf u)_t\psi \rd x+\int_\Omega \nabla\theta(\mathbf u)\nabla\psi\rd x=\mu \int_\Omega\theta(\mathbf u)(\nabla\cdot \mathbf u_t)\psi\rd x,~~~\forall~\psi\in H^1(\Omega).
\end{equation}
By taking $\psi=\theta(\mathbf u)$, we obtain
\begin{align*}
	\frac 12 \frac{\rd}{\rd t}\|\theta(\mathbf u)\|_{L^2(\Omega)}^2+\|\nabla\theta(\mathbf u)\|_{L^2(\Omega)}^2
	\leqslant|\mu|\cdot\|\nabla\cdot\mathbf{u}_t\|_{L^\infty(\Omega)}\|\theta(\mathbf u)\|_{L^2(\Omega)}^2.
\end{align*}
Notice that $V_n$ is a finite dimensional space, thus all norms are equivalent. It yields
\begin{equation}\label{Vnormequ}
	\|\nabla\cdot \mathbf u_t\|_{L^\infty(\Omega)}\leqslant C_n\|\mathbf u_t\|_{L^2(\Omega)}.
\end{equation}
With the help of the above inequality, we get
\begin{align}\label{theatnG}
	\frac 12 \frac{\rd}{\rd t}\|\theta(\mathbf u)\|_{L^2(\Omega)}^2+\|\nabla\theta(\mathbf u)\|_{L^2(\Omega)}^2
	\leqslant C_n(\|\mathbf u\|_{X_n})\|\theta(\mathbf u)\|_{L^2(\Omega)}^2.
\end{align}
By using the Gronwall's inequality, we have
\begin{align}\label{thetanL2}
	\|\theta(\mathbf u)\|_{L^2(\Omega)}^2\leqslant C_n(\|\mathbf u\|_{X_n})\|\theta_0\|_{L^2(\Omega)}^2,
\end{align}
 where $C_n(\|\mathbf u\|_{X_n})$ depends on $n$, but $n$ is fixed in this proof, which means this dependence does not influence the proof. Further, by integrating over $(0,t)$ in \eqref{theatnG}, we arrive at
\begin{align}\label{thetaH1bound}
	\theta(\mathbf u)\in L^2\left(0,T; H^1(\Omega)\right).
\end{align}
Now, with given $\theta(\mathbf u)$, we define $\mathcal{A}(\mathbf u)$ as a solution to the following initial boundary value problem
\begin{equation*}
	\begin{cases}
		(\mathcal{A}\mathbf u)_{tt}-\Delta(\mathcal{A}\mathbf u)-\Delta(\mathcal{A}\mathbf u)_{tt}=\mu\nabla \theta(\mathbf u),~~~~&\text{in}~\Omega\times(0,T),\\
		\mathcal{A}\mathbf u=0,&\text{on}~\partial\Omega\times(0,T),\\
		(\mathcal{A}\mathbf u)(0,\cdot )=u_0,~(\mathcal{A}\mathbf u)_t(0,\cdot )=v_0.
	\end{cases}
\end{equation*}	
More precisely, $\mathcal{A}\mathbf u$ satisfies the following equation
\begin{align}\label{Aulinear}
	\int_\Omega(\mathcal{A}\mathbf u)_{tt}\varphi\rd x+\int_\Omega\nabla(\mathcal{A}\mathbf u)\cdot\nabla\varphi\rd x+\int_\Omega\nabla(\mathcal{A}\mathbf u)_{tt}\cdot\nabla\varphi\rd x
	=-\mu\int_\Omega\theta(\mathbf u)\nabla\cdot\varphi\rd x,
\end{align}	
for all $\varphi\in V_n$.
The existence of $\mathcal{A}\mathbf u\in X_n$ is established by solving the linear ODE system \eqref{Aulinear} in the standard way. Consequently, the mapping $\mathcal{A}: X_n\rightarrow X_n$ is well-defined.

\textbf{Step 2.~We claim the mapping $\mathcal{A}$ is continuous and compact.
First, we prove $\mathcal{A}$ is continuous.}
We just need to prove that $\mathcal{A}$ can map convergent sequences to convergent sequences on $X_n$.
 Let us take $\{\mathbf u^m\}\in X_n$ be a convergent sequence and satisfy $\lim\limits_{m\rightarrow\infty}\mathbf u^m=\mathbf u$. Let $\sigma^m:=\theta(\mathbf u)-\theta(\mathbf u^m)$. By \eqref{theta-mathbf-u}, the following equality with zero initial data holds
\begin{align*}
	\frac{\rd}{\rd t}\int_\Omega\sigma^m\psi\rd x+\int_\Omega\nabla\sigma^m\cdot\nabla\psi\rd x
	=\mu\int_\Omega\theta(\mathbf u)\nabla\cdot\left(\mathbf u_t-\mathbf u_t^m\right)\psi\rd x+\mu \int_\Omega(\nabla\cdot \mathbf u_t^m)\sigma^m\psi\rd x.
\end{align*}
By taking $\psi=\sigma^m$, using \eqref{Vnormequ}, \eqref{thetanL2} and the Young inequality, we get
\begin{align*}
	\frac 12\frac{\rd}{\rd t}\int_\Omega|\sigma^m|^2\rd x+\int_\Omega|\nabla\sigma^m|^2\rd x
	&\leqslant|\mu|\cdot\|\nabla\cdot\left(\mathbf u_t-\mathbf u_t^m\right)\|_{L^\infty(\Omega)}\|\theta(\mathbf u)\|_{L^2(\Omega)}\|\sigma^m\|_{L^2(\Omega)}\\
	&~~~+|\mu|\cdot\|\nabla\cdot \mathbf u_t^m\|_{L^\infty(\Omega)}\|\sigma^m\|_{L^2(\Omega)}^2\\
	&\leqslant C_n(\|\mathbf u\|_{X_n})\|\nabla\cdot\left(\mathbf u_t-\mathbf u_t^m\right)\|_{L^\infty(\Omega)}^2+C_n(\|\mathbf u\|_{X_n})\|\sigma^m\|_{L^2(\Omega)}^2.
\end{align*}
Together with the Gronwall's inequality, we obtain
\begin{align}\label{sigmainfty}
  \|\nabla\sigma^m\|_{L^\infty\left((0,T); L^2(\Omega)\right)}^2
  \leqslant C_n(\|\mathbf u\|_{X_n})\|\nabla\cdot\left(\mathbf u_t-\mathbf u_t^m\right)\|_{L^\infty((0,T)\times\Omega)}^2,
\end{align}
where $C_n(\|\mathbf u\|_{X_n})>0$ is a constant depending on $X_n$, but independent of $m$.
Denoting $\phi^m:=\mathcal{A}\mathbf u-\mathcal A\mathbf u^m$ and using \eqref{Aulinear}, then we obtain that $\phi^m$ satisfies the following equation with zero initial data
\begin{align*}
		\int_\Omega(\phi^m)_{tt}\varphi\rd x+\int_\Omega\nabla(\phi^m)\cdot\nabla\varphi\rd x+\int_\Omega\nabla(\phi^m)_{tt}\cdot\nabla\varphi\rd x
		=-\mu\int_\Omega\sigma^m\nabla\cdot\varphi\rd x,
\end{align*}
for all $\varphi\in V_n$. By taking $\varphi=\phi_t^m$ in the above equality and using the Young inequality, we obtain that
\begin{align*}
	\frac 12\frac{\rd }{\rd t}\bigg(\|\phi^m_t\|_{L^2(\Omega)}^2&+\|\nabla\phi^m\|_{L^2(\Omega)}^2+\|\nabla\phi^m_t\|_{L^2(\Omega)}^2\bigg)\\
	&\leqslant |\mu|\cdot\|\nabla\sigma^m\|_{L^2(\Omega)}\|\phi_t^m\|_{L^2(\Omega)}\\
	&\leqslant C_n(\|\mathbf u\|_{X_n})\|\nabla\cdot\left(\mathbf u_t-\mathbf u_t^m\right)\|_{L^\infty((0,T)\times\Omega)}\|\phi_t^m\|_{L^2(\Omega)}\\
	&\leqslant C_n(\|\mathbf u\|_{X_n})\|\nabla\cdot\left(\mathbf u_t-\mathbf u_t^m\right)\|_{L^\infty((0,T)\times\Omega)}^2+\|\phi_t^m\|_{L^2(\Omega)}^2.
\end{align*}
Applying the Gronwall's inequality again, together with \eqref{sigmainfty}, it yields
\begin{align*}
	\|\phi^m_t\|_{L^2(\Omega)}^2+\|\nabla\phi^m\|_{L^2(\Omega)}^2&+\|\nabla\phi^m_t\|_{L^2(\Omega)}^2\\
	&\leqslant C_n(\|\mathbf u\|_{X_n})\|\nabla\cdot\left(\mathbf u_t-\mathbf u_t^m\right)\|_{L^\infty((0,T)\times\Omega)}^2\rightarrow0,~~~m\rightarrow\infty,
\end{align*}
which implies
$$\|\phi^m\|_{X_n}=\|\mathcal{A}\mathbf u-\mathcal A\mathbf u^m\|_{X_n}\rightarrow0,~~~m\rightarrow\infty.$$
Thus, $\mathcal{A}$ is a continuous mapping on $X_n$.

\textbf{Next, we prove that $\mathcal{A}$ is compact.}
Suffice it to say that $\mathcal{A}$ can map bounded sets to compact sets on $X_n$.
Let $\mathbf u\in X_n$ and $\|\mathbf u\|_{X_n}\leqslant R$, where $R$ is a positive constant.
By taking $\varphi=(\mathcal{A}\mathbf u)_{tt}$ in \eqref{Aulinear} and using the Young inequality, we get
\begin{align*}
	\|(\mathcal{A}\mathbf u)_{tt}\|_{L^2(\Omega)}^2+\|(\nabla\mathcal{A}\mathbf u)_{tt}\|_{L^2(\Omega)}^2
	&=\int_\Omega\Delta(\mathcal{A}\mathbf u)\cdot (\mathcal{A}\mathbf u)_{tt}+\mu \int_\Omega\nabla\theta(\mathbf u)\cdot (\mathcal{A}\mathbf u)_{tt}\rd x\\
	&\leqslant \frac 12\|(\mathcal{A}\mathbf u)_{tt}\|_{L^2(\Omega)}^2+ C\left(\|\Delta (\mathcal{A}\mathbf u)\|_{L^2(\Omega)}^2+\|\nabla\theta(\mathbf u)\|_{L^2(\Omega)}^2\right).
\end{align*}
This allows us to obtain that
\begin{align}\label{Autt}
		\|(\mathcal{A}\mathbf u)_{tt}\|_{L^2(\Omega)}^2+\|(\nabla\mathcal{A}\mathbf u)_{tt}\|_{L^2(\Omega)}^2
		\leqslant C\left(\|\Delta (\mathcal{A}\mathbf u)\|_{L^2(\Omega)}^2+\|\nabla\theta(\mathbf u)\|_{L^2(\Omega)}^2\right)
		\leqslant C_n(R).
\end{align}
Here we have used equivalence of norms on $X_n$ and continuity of $\mathcal{A}$ to bound the first term, and \eqref{thetaH1bound} to bound the second term. Together with \eqref{Autt}, the fact that $X_n$ is finite dimensional space immediately shows compactness of the mapping $\mathcal{A}$ on $X_n$ as a consequence of the embedding $H^1(0,T)\hookrightarrow C[0,T]$.

\textbf{Step 3.~Boundedness of the set $Y$.}~Let us introduce a notation
 $$Y=\left\{\mathbf u\in X_n;~\mathbf u=\lambda\mathcal{A}\mathbf u,~\text{for some}~\lambda\in(0,1]\right\}.$$
Then $\left(\frac 1\lambda\mathbf u,\theta(\mathbf u)\right)$ is a solution to the approximate problem in the sense of Definition \ref{definitionApproximate}. By Lemma \ref{Lemappenergy}, for all $t\in[0,T]$, the following inequality
	\begin{align*}
	&~~~\frac{1}{2\lambda}\int_\Omega |\mathbf u_t|^2\rd x+\frac{1}{2\lambda}\int_\Omega |\nabla \mathbf u|^2\rd x+\frac{1}{2\lambda}\int_\Omega|\nabla \mathbf u_t|^2\rd x+\int_\Omega\theta(\mathbf u) \rd x\leqslant CE_0^n,
\end{align*}
is satisfied, which implies
$$\|\mathbf u\|_{X_n}\leqslant CE_0^n,~~~\mathbf u\in Y.$$

Then, there exists a fixed point $\mathbf u$ of the mapping $\mathcal{A}$ by Schaefer's fixed point theorem.
Due to Lemma \ref{theta+}, $\theta(\mathbf u)$ is a strictly positive function.
Thus, $(\mathbf u,\theta(\mathbf u))$ is a solution of the approximate problem in the sense of Definition \ref{definitionApproximate}.
Finally, since $\theta\nabla\cdot\mathbf{u}\in L^2\left(0,T; L^2(\Omega)\right)$, we get $\theta(\mathbf u)\in L^2\left(0,T; H^2(\Omega)\right)\cap H^1\left(0,T; L^2(\Omega)\right)$ by the standard parabolic theory.
Thus, the proof of Proposition \ref{thetangeq0} is now finished.
\end{proof}
\begin{remark}\label{Re-Appregularity}
	Notice that $u^n\in V_n\subset C^\infty(\Omega)$. Thus, we get $u^n\in C^\infty(\Omega)$ for any $t\geqslant0$. Due to the regularity theory for Galerkin method \cite{Evans1}, we obtain $u^n\in C^{3,\infty}\left([0,T]\times \Omega\right)$. Similarly, we have $\theta^n\in C^{1,2}\left([0,T]\times \Omega\right)$ by standard regularity theory of parabolic equations \cite{FA1964,LSU1968}.
\end{remark}

\subsection{Passing to the limit}
The key observation is that the solutions to the approximate problem are smooth by Remark \ref{Re-Appregularity} and $\theta^n\geqslant c_n>0,~n\in\mathbb{N}$ (by Lemma \ref{theta+}), so we can replace $\psi$ by $(\theta^n)^{-1}\psi$ in the heat equation of \eqref{equAppro} and denote $\tau^n:=\log\theta^n$. Consequently, we can rewrite $\eqref{equAppro}_2$ in the form of entropy equation
\begin{equation}\label{equ-entropy}
	\int_\Omega\tau^n_t\psi\rd x+\int_\Omega\nabla\tau^n\nabla\psi\rd x-\int_\Omega(\nabla\tau^n)^2\psi\rd x=-\mu\int_\Omega u_t^n\cdot\nabla\psi\rd x,
\end{equation}
for all $\psi\in H^1(\Omega)$. Notice that the above entropy equation is equivalent to the heat equation in \eqref{equAppro} provided that $\theta^n>0$.

Next, we will prove a quantitative formulation of the second law of thermodynamics. This following result will be especially useful in deriving compactness estimates and analyzing the asymptotic behavior of solutions.

\begin{proposition}\label{ProtauC}
	Let $(u^n,\theta^n)$ be a weak solution of the problem \eqref{equAppro} in the sense of Definition \ref{definitionApproximate}. Then,  for almost all $t\in[0,\infty)$, the following equation
	\begin{align*}
		\frac{\rd}{\rd t}\int_\Omega\tau^n\rd x=\int_\Omega|\nabla\tau^n|^2\rd x,
	\end{align*}
	is satisfied, where $\tau^n:=\log\theta^n$.
	Moreover, there exists $C>0$ such that
	\begin{align}\label{nablatau}
		\int_0^\infty\int_\Omega|\nabla\tau^n|^2\rd x\rd t\leqslant C,
	\end{align}
	where the constant $C$ is positive and independent of time and $n$.
\end{proposition}
\begin{proof}
By Remark \ref{Re-Appregularity} the entropy equation \eqref{equ-entropy} is satisfied.
Thus, by taking $\psi=1$ in \eqref{equ-entropy}, we obtain
$$\frac{\rd}{\rd t}\int_\Omega\tau^n\rd x=\int_\Omega|\nabla\tau^n|^2\rd x.$$
Integrating over $(0,t)$ for an arbitrary $t>0$, we obtain
\begin{align}\label{tau0}
	\int_\Omega\tau^n\rd x
	=\int_0^t\int_\Omega|\nabla\tau^n|^2\rd x\rd s+\int_\Omega\tau_0^n\rd x,
\end{align}
where $\tau_0^n=\log\theta^n_0$. Combining \eqref{tau0} and Lemma \ref{Lemappenergy}, we arrive at
\begin{align*}
	\frac 12\int_\Omega |u_t^n|^2\rd x&+\frac 12\int_\Omega |\nabla u^n|^2\rd x+\frac 12\int_\Omega |\nabla u_t^n|^2\rd x+\int_0^t\int_\Omega|\nabla\tau^n|^2\rd x\rd s+\int_\Omega\theta^n \rd x-\int_\Omega\tau^n\rd x\\
	&=\frac 12\int_\Omega |v^n_0|^2\rd x+\frac 12\int_\Omega |\nabla u^n_{0}|^2\rd x+\frac 12\int_\Omega |\nabla v^n_0|^2\rd x+\int_\Omega\theta^n_0 \rd x-\int_\Omega\tau^n_0\rd x<\infty,
\end{align*}
since
$$\int_\Omega\theta^n \rd x>\int_\Omega\tau^n\rd x.$$
The claim is shown.	
\end{proof}
\begin{remark}
	We observe that the proofs of Proposition \ref{ProtauC} and Lemma \ref{Lemappenergy} imply that there exists a constant $C>0$, independent of both time and $n$, such that
	\begin{align}\label{tauL1}
		\int_\Omega|\tau^n(t,x)|\rd x< C,
	\end{align}
	which plays a crucial role in last section when we examine the convergence of $\theta$ as time approaches infinity.
\end{remark}

In order to complete the proof of Theorem \ref{ThexistandPos},  we still need to build up some useful lemmas.
\begin{lemma}\label{Ththetainfty}
	Let $\theta^n$ be a solution of the problem \eqref{equAppro} in the sense of Definition \ref{definitionApproximate}. Then we have
	$$\theta^n\in L^\infty\left((0,\infty)\times\Omega\right).$$
\end{lemma}
\begin{proof} The proof splits into some steps. 	
\textbf{Step 1. We claim that }
$$\theta^n\in L^\infty\left((0,\infty); L^2(\Omega)\right).$$
This step seems redundant since all the reasoning is contained in the iteration step. However, it is illustrative so that the readers see precisely what comes next. Let us take $\psi=\theta$ in the second equation of \eqref{equAppro} and use the H\"older inequality. Then, we get
\begin{align}\label{thetaL2}
	\frac 12\frac{\rd}{\rd t}\int_\Omega|\theta^n|^2\rd x+\int_\Omega|\nabla\theta^n|^2\rd x
	&=-2\mu\int_\Omega  u_t^n \theta^n\nabla\theta^n\rd x\notag\\
	&\leqslant2|\mu|\cdot\|u_t^n\|_{L^6(\Omega)}\|\theta^n\|_{L^3(\Omega)}\|\nabla\theta^n\|_{L^2(\Omega)}.
\end{align}
The norm $\|\theta^n\|_{L^3(\Omega)}$ is bounded via the Gagliardo-Nirenberg inequality
$$\|\theta^n\|_{L^3(\Omega)}\leqslant C\|\nabla\theta^n\|_{L^2(\Omega)}^{\frac 45}\|\theta^n\|_{L^1(\Omega)}^{\frac 15}+C\|\theta^n\|_{L^1(\Omega)}.$$
With the help of the above inequality, Sobolev embedding $H^1\hookrightarrow L^6$, the fact that $\|\theta^n\|_{L^1(\Omega)}$ and $\|\nabla u_t^n\|_{L^2(\Omega)}$ is bounded and independent of both time and $n$ from Lemma \ref{Lemappenergy}, and the Young inequality, we obtain that
\begin{align*}
	2|\mu|\cdot\|u_t^n\|_{L^6(\Omega)}&\|\theta^n\|_{L^3(\Omega)}\|\nabla\theta^n\|_{L^2(\Omega)}\\
	&\leqslant C\|u_t^n\|_{L^6(\Omega)}\cdot\left(\|\nabla\theta^n\|_{L^2(\Omega)}^{\frac 45}\|\theta^n\|_{L^1(\Omega)}^{\frac 15}+\|\theta^n\|_{L^1(\Omega)}\right)\|\nabla\theta^n\|_{L^2(\Omega)}\\
	&\leqslant C\|\nabla u_t^n\|_{L^2(\Omega)}\|\nabla\theta^n\|_{L^2(\Omega)}^{\frac 95}\|\theta^n\|_{L^1(\Omega)}^{\frac 15}+C\|\nabla u_t^n\|_{L^2(\Omega)}\|\theta^n\|_{L^1(\Omega)}\|\nabla\theta^n\|_{L^2(\Omega)}\\
	&\leqslant \frac 12\|\nabla\theta^n\|_{L^2(\Omega)}^2 +C\|\theta^n\|_{L^1(\Omega)}^2.
\end{align*}
Inserting the above inequality into \eqref{thetaL2}, it gives us that
\begin{align*}
	\frac{\rd}{\rd t}\int_\Omega|\theta^n|^2\rd x+\int_\Omega|\nabla\theta^n|^2\rd x\leqslant C.
\end{align*}
The Gagliardo-Nirenberg inequality, the Young inequality and Lemma \ref{Lemappenergy} ensure that
\begin{align*}
	\|\theta^n\|_{L^2(\Omega)}^2&\leqslant C\|\nabla\theta^n\|_{L^2(\Omega)}^{\frac 65}\|\theta^n\|_{L^1(\Omega)}^{\frac 45}+C\|\theta^n\|_{L^1(\Omega)}^2\\
	&\leqslant \|\nabla\theta^n\|_{L^2(\Omega)}^2+C\|\theta^n\|_{L^1(\Omega)}^2
	\leqslant \|\nabla\theta^n\|_{L^2(\Omega)}^2+C.
\end{align*}
Hence we get
\begin{align*}
	\frac{\rd}{\rd t}\int_\Omega|\theta^n|^2\rd x+\int_\Omega|\theta^n|^2\rd x\leqslant C.
\end{align*}
It immediately yields
\begin{align*}
	\int_\Omega|\theta^n|^2\rd x\leqslant \int_\Omega|\theta_0|^2\rd x+C.
\end{align*}
Since $C$ is independent of both time and $n$ in the above inequality, we obtain $$\theta^n\in L^\infty\left(0,\infty; L^2(\Omega)\right).$$

\textbf{Step 2. Estimate the norm $\|\theta\|_{L^\infty\left(0,\infty; L^{2^{m+1}}(\Omega)\right)},~m\in\mathbb{N}$.}
By taking $\psi=\theta^{2^{{m+1}}-1}$ in the second equation of \eqref{equAppro} and integrating over $\Omega$ by parts, we deduce that
\begin{align}\label{theta1}
	\frac{1}{2^{m+1}}\frac{\rd}{\rd t}\int_\Omega(\theta^n)^{2^{m+1}}\rd x+\frac{2^{m+1}-1}{4^m}\int_\Omega\left[\nabla\left((\theta^n)^{2^m}\right)\right]^2\rd x=\mu\int_\Omega\dive(u_t^n)\cdot(\theta^n)^{2^{{m+1}}}\rd x.
\end{align}
We shall estimate the right-hand side of the above equality. After integration by parts, Lemma \ref{Lemappenergy} with Sobolev embedding $H^1\hookrightarrow L^6$ ensures that
\begin{align}\label{thetanR}
	\mu\int_\Omega\dive(u_t^n)\cdot(\theta^n)^{2^{{m+1}}}\rd x
	&=-2\mu\int_\Omega u_t^n(\theta^n)^{2^m}\nabla\left((\theta^n)^{2^m}\right)\rd x\notag\\
	&\leqslant C\|u_t^n\|_{L^6(\Omega)}\|(\theta^n)^{2^m}\|_{L^3(\Omega)}\|\nabla\left((\theta^n)^{2^m}\right)\|_{L^2(\Omega)}\notag\\
	&\leqslant C\|\nabla u_t^n\|_{L^2(\Omega)}\|(\theta^n)^{2^m}\|_{L^3(\Omega)}\|\nabla\left((\theta^n)^{2^m}\right)\|_{L^2(\Omega)}\notag\\
	&\leqslant C\|(\theta^n)^{2^m}\|_{L^3(\Omega)}\|\nabla\left((\theta ^n)^{2^m}\right)\|_{L^2(\Omega)}.
\end{align}
The Gagliardo-Nirenberg inequality gives us that
\begin{align}\label{thetan3}
	\|(\theta^n)^{2^m}\|_{L^3(\Omega)}\leqslant C\|\nabla\left((\theta^n)^{2^m}\right)\|_{L^2(\Omega)}^{\frac 45}\|(\theta^n)^{2^m}\|_{L^1(\Omega)}^{\frac 15}+C\|(\theta^n)^{2^m}\|_{L^1(\Omega)}.
\end{align}
Plugging \eqref{thetan3} into \eqref{thetanR} and utilizing the Young inequality, it yields
\begin{align*}
	\mu\int_\Omega\dive&(u_t^n)\cdot(\theta^n)^{2^{{m+1}}}\rd x\\
	&\leqslant  C\left(\|\nabla\left((\theta^n)^{2^m}\right)\|_{L^2(\Omega)}^{\frac 45}\|(\theta^n)^{2^m}\|_{L^1(\Omega)}^{\frac 15}+\|(\theta^n)^{2^m}\|_{L^1(\Omega)}\right)\|\nabla\left((\theta ^n)^{2^m}\right)\|_{L^2(\Omega)}\\
	&\leqslant C\|\nabla\left((\theta^n)^{2^m}\right)\|_{L^2(\Omega)}^{\frac 95}\|(\theta^n)^{2^m}\|_{L^1(\Omega)}^{\frac 15}+C\|(\theta^n)^{2^m}\|_{L^1(\Omega)}\|\nabla\left((\theta ^n)^{2^m}\right)\|_{L^2(\Omega)}\\
	&\leqslant \varepsilon\|\nabla\left((\theta^n)^{2^m}\right)\|_{L^2(\Omega)}^2+C\left(\varepsilon^{-1}+\varepsilon^{-9}\right)\|(\theta^n)^{2^m}\|_{L^1(\Omega)}^{2},
\end{align*}
where $\varepsilon>0$ will be defined later. Inserting the above inequality into \eqref{theta1}, we obtain
\begin{align*}
	\frac{1}{2^{m+1}}\frac{\rd}{\rd t}\int_\Omega(\theta^n)^{2^{m+1}}\rd x
	&+\alpha_m\int_\Omega\left[\nabla\left((\theta^n)^{2^m}\right)\right]^2\rd x\\
	&\leqslant \varepsilon\|\nabla\left((\theta^n)^{2^m}\right)\|_{L^2(\Omega)}^2+C\left(\varepsilon^{-1}+\varepsilon^{-9}\right)\|(\theta^n)^{2^m}\|_{L^1(\Omega)}^{2}.
\end{align*}
Here we denoted $\alpha_m:=\frac{2^{m+1}-1}{4^m}$.
Choosing $\varepsilon=\frac{\alpha_m}{2}$ in the above inequality, we get
\begin{align*}
\frac{1}{2^{m+1}}\frac{\rd}{\rd t}\int_\Omega(\theta^n)^{2^{m+1}}\rd x
+\frac{\alpha_m}{2}\int_\Omega\left[\nabla\left((\theta^n)^{2^m}\right)\right]^2\rd x
   \leqslant C\left(\alpha_m^{-1}+\alpha_m^{-9}\right)\|(\theta^n)^{2^m}\|_{L^1(\Omega)}^{2}.
\end{align*}
From the Gagliardo-Nirenberg inequality and the Young inequality, we arrive at
\begin{align*}
	\frac{1}{2^{m+1}}\frac{\rd}{\rd t}\int_\Omega(\theta^n)^{2^{m+1}}\rd x
	+\frac{\alpha_m}{2}\int_\Omega\left[(\theta^n)^{2^m}\right]^2\rd x
	\leqslant C\left(\alpha_m^{-1}+\alpha_m^{-9}+1\right)\|(\theta^n)^{2^m}\|_{L^1(\Omega)}^{2}.
\end{align*}
Since $\frac{\alpha_m}{2}\cdot2^{m+1}=2-\frac{1}{2^{m}}>1$, we can rewrite the above inequality as follows
\begin{align*}
	\frac{\rd}{\rd t}\|\theta^n\|_{L^{2^{m+1}}(\Omega)}^{2^{m+1}}
	+\|\theta^n\|_{L^{2^{m+1}}(\Omega)}^{2^{m+1}}
	\leqslant C_1\|\theta^n\|_{L^{2^{m}}(\Omega)}^{2^{m+1}},
\end{align*}
where $C\left(\alpha_m^{-1}+\alpha_m^{-9}+1\right)2^{m+1}\leqslant C2^{10m}=C_1$. 
Then we have
\begin{align*}
	\|\theta^n\|_{L^{2^{m+1}}(\Omega)}^{2^{m+1}}
	&\leqslant\re^{-t}\|\theta_0^n\|_{L^{2^{m+1}}(\Omega)}^{2^{m+1}}
	+C_1\|\theta^n\|_{L^{2^{m}}(\Omega)}^{2^{m+1}}\\
	&\leqslant|\Omega|\cdot\|\theta^n_0\|_{L^{\infty}(\Omega)}^{2^{m+1}}+C2^{10m}\|\theta^n\|_{L^\infty\left(0,\infty; L^{2^{m}}(\Omega)\right)}^{2^{m+1}},
\end{align*}
which indicates that
\begin{align*}
	\|\theta^n\|_{L^\infty\left(0,\infty; L^{2^{m+1}}(\Omega)\right)}
	\leqslant\left(|\Omega|\cdot\|\theta_0^n\|_{L^{\infty}(\Omega)}^{2^{m+1}}+C2^{10m}\|\theta^n\|_{L^\infty\left(0,\infty; L^{2^{m}}(\Omega)\right)}^{2^{m+1}}+1\right)^\frac{1}{2^{m+1}}.
\end{align*}
	
\textbf{Step 3. Start iteration.} We denote
$$\beta_m=\max\left\{\|\theta^n_0\|_{L^\infty(\Omega)},~\|\theta^n\|_{L^{\infty}\left(0,\infty; L^{2^m}(\Omega)\right)},~1\right\}.$$
Then, we obtain
$$\beta_{m+1}\leqslant \beta_mC^{\frac{1}{2^{m+1}}}2^\frac{10m}{2^{m+1}},$$
which implies that
$$\beta_m\leqslant \beta_1\prod_{k=2}^{m}C^{\frac{1}{2^{k}}}2^\frac{10(k-1)}{2^{k}}=\beta_1C^{\sum_{k=2}^{m}{\frac{1}{2^{k}}}}2^{\sum_{k=2}^{m}\frac{10(k-1)}{2^{k}}}
\leqslant\beta_1C^{\sum_{k=2}^{\infty}{\frac{1}{2^{k}}}}2^{\sum_{k=2}^{\infty}\frac{10(k-1)}{2^{k}}}.$$
Convergence of the two series in the exponents implies uniform boundedness $\beta_m$ for all $m\in\mathbb{N}$. Therefore, we get
$$\|\theta^n\|_{L^\infty\left(0,\infty; L^{2^m}(\Omega)\right)}\leqslant C,~~~\text{for all}~m\in\mathbb{N}.$$
The proof of Lemma \ref{Ththetainfty} is now finished.
\end{proof}

Lemma \ref{Ththetainfty} and Proposition \ref{ProtauC} give us immediately the following result, which plays a key role to estimate higher-order derivatives of $u^n$ and $\theta^n$. This following result provides a crucial technical tool to handle a tricky Gronwall-type inequality, ultimately yielding the required compactness estimates.
\begin{corollary}
	Let $\theta^n$ be a positive weak solution of \eqref{equAppro} in the sense of Definition \ref{definitionApproximate}. Then, we have
	\begin{align}\label{sqrttheta}
		\int_0^\infty\int_\Omega\frac{|\nabla\theta^n|^2}{\theta^n}\rd x\rd t<\infty.
	\end{align}
\end{corollary}
\begin{proof}
	From \eqref{nablatau} and Lemma \ref{Ththetainfty}, we obtain
	\begin{align*}
		\int_0^\infty\int_\Omega\frac{|\nabla\theta^n|^2}{\theta^n}\rd x\rd t
		\leqslant\|\theta^n\|_{L^\infty((0,\infty)\times\Omega)}\int_0^\infty\int_\Omega|\nabla\tau^n|^2\rd x\rd t\leqslant C.
	\end{align*}
	Thus, \eqref{sqrttheta} follows.
\end{proof}
Next, we present the key identity, which leads to the global existence of solutions. The approach parallels the one-dimensional treatment \cite{CTglobal1D}. For the reader's convenience, we present the details.
\begin{lemma} \label{Lemhigerpre}
Assume that $u^n$ and $\theta^n,~\theta^n>0$ are regular enough solutions of \eqref{equAppro} given by Definition \ref{definitionApproximate}. Then the following identity holds
	\begin{align}\label{higher1}
		\frac{1}{2}\frac{\rd}{\rd t}\bigg(\int_\Omega \frac{|\nabla \theta^n|^2}{\theta^n}\rd x&+\int_\Omega |\dive (u_t^n)|^2\rd x
		+\int_\Omega |\nabla \dive(u^n)|^2\rd x+\int_\Omega|\nabla\dive(u_t^n)|^2\rd x\bigg)\notag\\
		&=-\int_\Omega \theta^n \left|D^2\log\theta^n\right|^2\rd x
		+\frac{\mu}{2}\int_\Omega \frac{|\nabla \theta^n|^2}{\theta^n}\dive(u_t^n)\rd x.
	\end{align}
\end{lemma}
\begin{proof}
We omit the superscript for $\theta^n$ and $u^n$ in this proof.
By replacing $\psi$ by $(\sqrt{\theta})^{-1}\psi$ in the heat equation of \eqref{equAppro}, we obtain that
\begin{align*}
	\int_\Omega\frac{\theta_t}{\sqrt{\theta}}\psi\rd x-\int_\Omega\frac{\Delta\theta}{\sqrt{\theta}}\psi\rd x=\mu\int_\Omega\sqrt{\theta}\dive(u_t)\psi\rd x.
\end{align*}
Then, a simple computation shows that
\begin{align}\label{sqrt1}
	\frac{\rd}{\rd t}\int_\Omega\sqrt{\theta}\psi\rd x
	=\int_\Omega\Delta\sqrt{\theta}\cdot\psi\rd x+\int_\Omega\frac{|\nabla\sqrt\theta|^2}{\sqrt\theta}\psi\rd x
	+\frac \mu2\int_\Omega\sqrt{\theta}\dive(u_t)\psi\rd x.
\end{align}
Here we have used the following equality
\begin{align}\label{D2}
		|\nabla\theta|^2=4\theta|\nabla\sqrt\theta|^2,
	~~~\Delta\sqrt{\theta}=\nabla\cdot\left(\frac{\nabla\theta}{2\sqrt{\theta}}\right)=\frac{\Delta\theta}{2\sqrt{\theta}}-\frac{|\nabla\theta|^2}{4\theta^\frac32}.
\end{align}
Let us take $\psi=\Delta\sqrt\theta$ in \eqref{sqrt1}.
Together with \eqref{D2}, we get
\begin{align*}
	-\frac 12&\frac{\rd}{\rd t}\int_\Omega|\nabla\sqrt\theta|^2\rd x\\
	&=\int_\Omega|\Delta\sqrt\theta|^2\rd x-\int_\Omega\nabla\left(\frac{|\nabla\sqrt\theta|^2}{\sqrt\theta}\right)\nabla\sqrt\theta\rd x+\frac \mu2\int_\Omega\sqrt\theta\dive(u_t)\Delta\sqrt\theta\rd x\notag\\
	&=\int_\Omega|D^2\sqrt\theta|^2\rd x-2\int_\Omega\frac{|\nabla\sqrt\theta|^2D^2\sqrt\theta}{\sqrt\theta}\rd x+\int_\Omega\frac{|\nabla\sqrt\theta|^4}{\theta}\rd x
	+\frac \mu2\int_\Omega\sqrt\theta\dive(u_t)\Delta\sqrt\theta\rd x\notag\\
	&=\int_\Omega\left(D^2\sqrt\theta-\frac{|\nabla\sqrt\theta|^2}{\sqrt\theta}\right)^2\rd x+\frac \mu4\int_\Omega\Delta\theta\dive(u_t)\rd x-\frac \mu8\int_\Omega\frac{|\nabla\theta|^2}{\theta}\dive(u_t)\rd x\notag\\
	&=\frac 14\int_\Omega\theta\left(D^2\log\theta\right)^2\rd x-\frac \mu4\int_\Omega\nabla\theta\nabla\dive(u_t)\rd x-\frac \mu8\int_\Omega\frac{|\nabla\theta|^2}{\theta}\dive(u_t)\rd x.
\end{align*}
	Using \eqref{D2} again, we obtain
	\begin{align}\label{sqrt2}
		\frac 12\frac{\rd}{\rd t}\int_\Omega\frac{|\nabla\theta|^2}{\theta}\rd x
		&=-\int_\Omega\theta\left(D^2\log\theta\right)^2\rd x+\frac \mu2\int_\Omega\frac{|\nabla\theta|^2}{\theta}\dive(u_t)\rd x\notag\\
		&~~~+\mu\int_\Omega\nabla\theta\nabla\dive(u_t)\rd x.
	\end{align}
	Next, we replace $\varphi=\nabla\dive(u_t)$ in the first equation of \eqref{equAppro} and integrate over $\Omega$. It yields
	\begin{align}\label{diveut1}
	\frac 12\frac{\rd}{\rd t}\bigg(\int_\Omega |\dive(u_t)|^2\rd x+\int_\Omega|\nabla\dive(u)|^2\rd x&+\int_\Omega|\nabla\dive(u_t)|^2\rd x\bigg)\notag\\
		&=-\mu\int_\Omega\nabla\theta\nabla\dive(u_t)\rd x.
	\end{align}
By adding \eqref{sqrt2} and \eqref{diveut1}, the claim follows.	
\end{proof}

In order to obtain the limit of the approximate problem, compactness estimates are required. In what follows, we establish the estimates for higher-order derivatives of both displacement and temperature in the following theorem.
\begin{theorem}\label{Thhigher}
Let $(u^n,\theta^n)$ be the solution of the approximate problem in the sense of Definition \ref{definitionApproximate}.
Then, the following estimates (independent of both time and $n$) are satisfied
\begin{align*}
	\int_\Omega |\nabla \theta^n|^2\rd x+\int_\Omega |\dive (u_t^n)|^2\rd x+\int_\Omega |\nabla \dive(u^n)|^2\rd x+\int_\Omega|\nabla\dive(u_t^n)|^2\rd x\leqslant C,
\end{align*}
where the constant $C$ only depends on $\mu,~\theta_0,~\Omega,~\|\theta_0\|_{H^1(\Omega)},~\|u_0\|_{H^2(\Omega)},\|v_0\|_{H_0^1(\Omega)},~\|v_0\|_{H^2(\Omega)}$.
\end{theorem}
\begin{proof}
For the sake of simplicity, we omit the superscripts for $u^n$ and $\theta^n$ in the proof. We shall estimate the last term on the right-hand side of \eqref{higher1}. By using \eqref{D2}, Sobolev embedding $H^1\hookrightarrow L^6$ and the H\"older inequality, we obtain that
\begin{align}\label{higher2}
	\frac{\mu}{2}\int_\Omega \frac{|\nabla \theta |^2}{\theta}\dive(u_t)\rd x
	&\leqslant 2|\mu|\int_\Omega |\nabla\sqrt\theta |^2\dive(u_t)\rd x\notag\\
	&\leqslant C\|\dive(u_t)\|_{L^6(\Omega)}\|\nabla\sqrt{\theta}\|_{L^3(\Omega)}\|\nabla\sqrt{\theta}\|_{L^2(\Omega)}\notag\\
	&\leqslant C\|\nabla\dive(u_t)\|_{L^2(\Omega)}\|\nabla\sqrt{\theta}\|_{L^3(\Omega)}\|\nabla\sqrt{\theta}\|_{L^2(\Omega)}.
\end{align}
Applying the Gagliardo-Nirenberg inequality to the norm $\|\nabla\sqrt{\theta}\|_{L^3(\Omega)}$, we get
\begin{align*}
	\|\nabla\sqrt{\theta}\|_{L^3(\Omega)}\leqslant	C\|D^2\sqrt{\theta}\|_{L^2(\Omega)}^{\frac 12}	\|\nabla\sqrt{\theta}\|_{L^2(\Omega)}^{\frac 12}+C\|\nabla\sqrt{\theta}\|_{L^2(\Omega)}.
\end{align*}
Plugging the above inequality into \eqref{higher2}, combining the Young inequality and Lemma \ref{LemDsqrt}, we deduce that
\begin{align}\label{higher3}
	\frac{\mu}{2}\int_\Omega &\frac{|\nabla \theta |^2}{\theta}\dive(u_t)\rd x\notag\\
	&\leqslant C\|\nabla\dive(u_t)\|_{L^2(\Omega)}\left(\|D^2\sqrt{\theta}\|_{L^2(\Omega)}^{\frac 12}	\|\nabla\sqrt{\theta}\|_{L^2(\Omega)}^{\frac 12}+\|\nabla\sqrt{\theta}\|_{L^2(\Omega)}\right)\|\nabla\sqrt{\theta}\|_{L^2(\Omega)}\notag\\
	&\leqslant C\|\nabla\dive(u_t)\|_{L^2(\Omega)}\|D^2\sqrt{\theta}\|_{L^2(\Omega)}^{\frac 12}	 \|\nabla\sqrt{\theta}\|_{L^2(\Omega)}^{\frac 32}+C\|\nabla\dive(u_t)\|_{L^2(\Omega)}\|\nabla\sqrt{\theta}\|_{L^2(\Omega)}^2\notag\\
	&\leqslant C\|D^2\sqrt{\theta}\|_{L^2(\Omega)}	\|\nabla\sqrt{\theta}\|_{L^2(\Omega)}+C\|\nabla\dive(u_t)\|_{L^2(\Omega)}^2 \|\nabla\sqrt{\theta}\|_{L^2(\Omega)}^2+C\|\nabla\sqrt{\theta}\|_{L^2(\Omega)}^2\notag\\
	&\leqslant\varepsilon\|D^2\sqrt{\theta}\|_{L^2(\Omega)}^2+C\|\nabla\sqrt{\theta}\|_{L^2(\Omega)}^2+C\|\nabla\dive(u_t)\|_{L^2(\Omega)}^2 \|\nabla\sqrt{\theta}\|_{L^2(\Omega)}^2\notag\\
	&\leqslant\varepsilon\cdot\frac{11+4\sqrt{3}}{8}\int_\Omega\theta[D^2(\log\theta)]^2\rd x+C_1\int_\Omega \frac{|\nabla \theta |^2}{\theta }\rd x\notag\\
	&~~~+C_2\int_\Omega|\nabla\dive(u_t)|^2\rd x\cdot\int_\Omega \frac{|\nabla \theta |^2}{\theta }\rd x.
\end{align}
Choosing $\varepsilon=\frac{4}{11+4\sqrt{3}}$ and substituting \eqref{higher3} into \eqref{higher1}, it yields
\begin{align}\label{higher4}
	\frac{\rd}{\rd t}\bigg(\int_\Omega \frac{|\nabla \theta |^2}{\theta }\rd x&+\int_\Omega |\dive (u_t)|^2\rd x+\int_\Omega |\nabla \dive(u)|^2\rd x+\int_\Omega|\nabla\dive(u_t)|^2\rd x\bigg)\notag\\
	&\leqslant C_1\int_\Omega \frac{|\nabla \theta |^2}{\theta }\rd x+C_2\int_\Omega|\nabla\dive(u_t)|^2\rd x\cdot\int_\Omega \frac{|\nabla \theta |^2}{\theta }\rd x.
\end{align}
Let us introduce the following notations
\begin{align*}
	&x(t)=\int_\Omega \frac{|\nabla \theta |^2}{\theta }\rd x,\\
	&y(t)=\int_\Omega|\nabla\dive(u_t)|^2\rd x,\\
	&z(t)=\int_\Omega |\dive (u_t)|^2\rd x+\int_\Omega |\nabla \dive(u)|^2\rd x.
\end{align*}
Then, we can rewrite \eqref{higher4} as follows
\begin{align*}
	\big(x(t)+y(t)+z(t)\big)'
	&\leqslant C_1x(t)+C_2x(t)y(t)\\
	&\leqslant C_1x(t)+C_2x(t)\big(x(t)+y(t)+z(t)\big).
\end{align*}
Multiplying both sides by $\re^{-C_2\int_0^sx(\tau)\rd\tau}$ and integrating over $(0,t)$ for arbitrary $t$, we get
\begin{align*}
	x(t)+y(t)+z(t)\leqslant
	\re^{C_2\int_0^tx(\tau)\rd\tau}\left(C_1\int_0^t\re^{-C_2\int_0^sx(\tau)\rd\tau}x(s)\rd s+x(0)+y(0)+z(0)\right).
\end{align*}
With the help of \eqref{sqrttheta}, we have
$$x(t)+y(t)+z(t)\leqslant C.$$	
In light of Lemma \ref{Ththetainfty} and \eqref{sqrttheta} again, we additionally obtain
\begin{align*}
	\|\nabla\theta\|_{L^2(\Omega)}^2
	\leqslant\|\theta\|_{L^{\infty}\left((0,\infty)\times\Omega\right)}\int_\Omega\frac{|\nabla\theta|^2}{\theta}\rd x
	\leqslant \|\theta\|_{L^{\infty}\left((0,\infty)\times\Omega\right)}x(t)\leqslant C,
\end{align*}
where $C$ is a positive constant independent of both the time and $n$. Thus, this completes the proof of Theorem \ref{Thhigher}.
\end{proof}

However, the weak star convergence results from Theorem \ref{Thhigher} cannot be directly applied to pass to the limit in \eqref{equAppro}, due to the presence of nonlinear terms. Consequently, the estimates below become indispensable.
\begin{theorem}\label{Thutt}
	Let $(u^n,\theta^n)$ be the solution of the approximate problem in the sense of Definition \ref{definitionApproximate}.
	Then, for all  $n\in\mathbb{N}$ we have
	$$u^n_{tt}\in L^\infty\left(0,\infty; H_0^1(\Omega)\right),~~~\theta_t^n\in L^2_{\rm loc}\left(0,\infty; L^2(\Omega)\right).$$
\end{theorem}
\begin{proof}
For the sake of simplicity, we omit the subscript for $\theta^n$ and $u^n$ in this proof.	
Taking $\psi=\theta_t$ in the heat equation of \eqref{equAppro} and integrating over $[0,T]$, we obtain
\begin{align}\label{highhert1}
	\int_0^T\int_\Omega|\theta_t|^2\rd x\rd t+\frac 12\int_\Omega|\nabla\theta|^2\rd x
	=\frac 12\int_\Omega|\nabla\theta_0|^2\rd x+\mu \int_0^T\int_\Omega\theta(\nabla\cdot u_t)\theta_t\rd x\rd t.
\end{align}
We shall calculate the right-hand side of the above equality. From Lemma \ref{Ththetainfty}, the H\"older inequality and the Young inequality, we conclude that
	\begin{align*}
		\mu \int_0^T\int_\Omega\theta(\nabla&\cdot u_t)\theta_t\rd x\rd t\\
		&\leqslant |\mu|\cdot\|\theta\|_{L^\infty\left((0,\infty)\times\Omega\right)}\left(\int_0^T\int_\Omega|\dive(u_t)|^2\rd x\rd t\right)^{\frac 12}
		\left(\int_0^T\int_\Omega|\theta_t|^2\rd x\rd t\right)^{\frac 12}\\
		&\leqslant\frac 12\int_0^T\int_\Omega|\theta_t|^2\rd x\rd t
		+C\int_0^T\int_\Omega|\dive(u_t)|^2\rd x\rd t.
	\end{align*}
Plugging the above inequality into \eqref{highhert1} and using Theorem \ref{Thhigher}, we obtain
\begin{align*}
	\int_0^T\int_\Omega|\theta_t|^2\rd x\rd t+\int_\Omega|\nabla\theta|^2\rd x
	\leqslant\int_\Omega|\nabla\theta_0|^2\rd x+C\int_0^T\int_\Omega|\dive(u_t)|^2\rd x\rd t,
\end{align*}
where $C$ does not depend on $n$. Finally, let us take $\varphi=u_{tt}$ in \eqref{equAppro} and use the H\"older inequality and the Young inequality. It yields that
	\begin{align*}
		\|u_{tt}\|_{L^2(\Omega)}^2+\|\nabla u_{tt}\|_{L^2(\Omega)}^2
		&=\int_\Omega\Delta u\cdot u_{tt}+\mu \int_\Omega\nabla\theta\cdot u_{tt}\rd x\\
		&\leqslant \frac 12\|u_{tt}\|_{L^2(\Omega)}^2+ C\left(\|\Delta u\|_{L^2(\Omega)}^2+\|\nabla\theta\|_{L^2(\Omega)}^2\right).
	\end{align*}
According to Theorem \ref{Thhigher}, we arrive at
	\begin{align*}
		\|u_{tt}\|_{L^2(\Omega)}^2+\|\nabla u_{tt}\|_{L^2(\Omega)}^2
		\leqslant C\left(\|\Delta u\|_{L^2(\Omega)}^2+\|\nabla\theta\|_{L^2(\Omega)}^2\right)\leqslant C,
	\end{align*}
where $C$ is a constant independent of both $n$ and the time variable. Thus, the proof of Theorem \ref{Thutt} now is completed.	
\end{proof}

We now prove that the approximation solution defined by Definition \ref{definitionApproximate} converges to the solution of \eqref{equ1.1} in the sense of Definition \ref{solution} based on the theorems above, which establishes the existence part of Theorem \ref{ThexistandPos}.
\begin{theorem}\label{Thexist1}
	For any $T>0$ and any constant $\mu\in\mathbb{R}$, there exists a global-in-time solution of \eqref{equ1.1} in the sense of Definition \ref{solution}.
\end{theorem}
\begin{proof}
By Theorem \ref{Thhigher} and Theorem \ref{Thutt}, there exist a subsequence of $(u^n,\theta^n)$ (still denoted as $(u^n,\theta^n)$) and
\begin{align*}
	&u\in L^\infty\left(0,\infty; H^2(\Omega)\right)\cap W^{1,\infty}\left(0,\infty; H_0^1(\Omega)\cap H^2(\Omega)\right)\cap W^{2,\infty}\left(0,\infty; H_0^1(\Omega)\right),\\
	&\theta\in L^\infty\left(0,\infty; H^1(\Omega)\right)\cap H^1_{\rm loc}\left(0,\infty; L^2(\Omega)\right),
\end{align*}
such that the following weak convergences
\begin{align}\label{weak*}
	\begin{split}
	u^n&\stackrel *{\rightharpoonup }u~~~~~~~\text{in}~L^{\infty }(0,\infty;H^2(\Omega)),\\
	u^n_t&\stackrel *{\rightharpoonup }u_t~~~~~~\text{in}~L^{\infty }(0,\infty;H^1_0(\Omega)\cap H^2(\Omega)),\\
	u^n_{tt}&\stackrel *{\rightharpoonup }u_{tt}~~~~~\text{in}~L^{\infty }(0,\infty;H_0^{1}(\Omega)),\\
	\theta^{n}&\stackrel *{\rightharpoonup }\theta~~~~~~~\text{in}~L^{\infty }(0,\infty;H^1(\Omega)),\\ \theta^n_{t}&\stackrel{}{\rightharpoonup }\theta_t~~~~~~\text{in}~L_{\rm loc}^2(0,\infty;L^{2}(\Omega)),
\end{split}
\end{align}
hold. By the Aubin-Lions Lemma, we obtain
\begin{align}\label{C0T}
	u\in C\left([0,\infty); H_0^1(\Omega)\right),~~~u_t\in C\left([0,\infty); H_0^1(\Omega)\right),~~~\theta\in C\left([0,\infty); L^2(\Omega)\right).
\end{align}
In the standard way, we know that $(u,\theta)$ satisfies the moment equation in Definition \ref{solution}.
Let us take $\psi\in C^\infty(\Omega)$, multiply the second equation in Definition \ref{definitionApproximate} and integrate over $[0,T]$. It yields
\begin{align*}
	\int_0^T\int_\Omega\theta^n_t\psi \rd x\rd t+\int_0^T\int_\Omega \nabla\theta^n\nabla\psi\rd x\rd t=\mu \int_0^T\int_\Omega\theta^n(\nabla\cdot u_t^n)\psi\rd x\rd t.
\end{align*}
From the convergences of \eqref{weak*} and \eqref{C0T}, it is easy to see that
\begin{align*}
	\int_0^T\int_\Omega\theta^n_t\psi \rd x\rd t+\int_0^T\int_\Omega \nabla\theta^n\nabla\psi\rd x\rd t
	\rightarrow \int_0^T\int_\Omega\theta_t\psi \rd x\rd t+\int_0^T\int_\Omega \nabla\theta\nabla\psi\rd x\rd t,
\end{align*}	
and
\begin{align*}
	\mu \int_0^T\int_\Omega\theta^n(\nabla\cdot u_t^n)\psi\rd x\rd t
	\rightarrow
	\mu \int_0^T\int_\Omega\theta(\nabla\cdot u_t)\psi\rd x\rd t,
\end{align*}	
since \eqref{weak*} and the Aubin-Lions Lemma give us that
\begin{align}\label{thetaL2C}
	\theta^n\rightarrow\theta~~~\text{in}~L^2\left(0,T; L^2(\Omega)\right).
\end{align}
This completes the proof of Theorem \ref{Thexist1}.	
\end{proof}

\subsection{Positivity of the temperature}
We claim that the temperature $\theta$ is positive in Theorem \ref{ThexistandPos}.
In what follows, we provide a complete proof of the claim. The proof of the lower bound relies on a Moser-type iteration technique on the negative part of the logarithm of temperature.
\begin{theorem}\label{Thtauinfty}
	Let $(u,\theta)$ be a solution by solving \eqref{equ1.1} in the sense of Definition \ref{solution}. Then, there exists a constant $\theta^*>0$ such that the temperature function satisfies
	$$\theta(t,x)\geqslant\theta^*>0,$$
	for almost all $(t,x)\in[0,\infty)\times\Omega$.
\end{theorem}
\begin{proof}
Let $\theta^n$ be a solution to the approximate problem in the sense of Definition \ref{definitionApproximate}. We still denote $\tau^n=\log\theta^n$. The proof splits into some steps.

\textbf{Step 1. We claim that $$\tau^n_-\in L^\infty\left((0,\infty); L^2(\Omega)\right),$$}
where $\tau^n_-:=\max\{0,-\tau^n\}$. Taking $\psi=-\tau^n_-$ in \eqref{equ-entropy}, we get
\begin{align*}
	\frac 12\frac{\rd}{\rd t}\int_\Omega|\tau^n_-|^2\rd x+\int_\Omega|\nabla(\tau^n_-)|^2\rd x+\int_\Omega|\nabla(\tau^n)|^2\tau^n_-\rd x
	=-\mu\int_\Omega\tau^n_-\cdot\dive(u_t^n)\rd x.
\end{align*}	
By using the H\"older inequality and the Young inequality, we obtain
\begin{align*}
	\frac 12\frac{\rd}{\rd t}\int_\Omega|\tau^n_-|^2\rd x+\int_\Omega|\nabla(\tau^n_-)|^2\rd x
	&\leqslant\mu\int_\Omega u_t^n\cdot\nabla(\tau^n_-)\rd x\\
	&\leqslant|\mu|\cdot\|u_t^n\|_{L^2(\Omega)}\|\|\nabla(\tau^n_-)\|_{L^2(\Omega)}\\
	&\leqslant\frac 12\|\nabla(\tau^n_-)\|_{L^2(\Omega)}^2+C\|u_t^n\|_{L^2(\Omega)}^2,
\end{align*}	
which implies	
\begin{align}\label{tau-1}
\frac{\rd}{\rd t}\int_\Omega|\tau^n_-|^2\rd x+\int_\Omega|\nabla(\tau^n_-)|^2\rd x
\leqslant C\|u_t^n\|_{L^2(\Omega)}^2.
\end{align}	
The Gagliardo-Nirenberg inequality, the Young inequality, Lemma \ref{Lemappenergy} and \eqref{tauL1} ensure that
\begin{align*}
	\|\tau^n_-\|_{L^2(\Omega)}^2&\leqslant C\|\nabla(\tau^n_-)\|_{L^2(\Omega)}^{\frac 65}\|\tau^n_-\|_{L^1(\Omega)}^{\frac 45}+C\|\tau^n_-\|_{L^1(\Omega)}^2\\
	&\leqslant \|\nabla(\tau^n_-)\|_{L^2(\Omega)}^2+C\|\tau^n_-\|_{L^1(\Omega)}^2
	\leqslant \|\nabla(\tau^n_-)\|_{L^2(\Omega)}^2+C.
\end{align*}	
Plugging the above inequality into \eqref{tau-1} and using Lemma \ref{Lemappenergy} again, it leads to
\begin{align*}
	\frac{\rd}{\rd t}\|\tau^n_-\|_{L^(\Omega)}^2+\|\tau^n_-\|_{L^(\Omega)}^2\leqslant C,
\end{align*}	
where $C$ is independent of both $n$ and time variable.
Thus, it immediately gives us
$$\tau^n_-\in L^\infty\left((0,\infty); L^2(\Omega)\right).$$
	
\textbf{Step 2. Estimate the norm $\|\tau^n_-\|_{L^\infty\left(0,\infty; L^{2^{m+1}}(\Omega)\right)},~m\in\mathbb{N}$.}
Taking $\psi=(\tau^n_-)^{2^{{m+1}-1}}$ and integrating over $\Omega$. By a straightforward calculation, we get
\begin{align}\label{tau-2}
	\frac{1}{2^{m+1}}\frac{\rd}{\rd t}\int_\Omega(\tau^n_-)^{2^{m+1}}\rd x+\frac{2^{m+1}-1}{4^m}\int_\Omega\left[\nabla\left((\tau^n_-)^{2^m}\right)\right]^2\rd x=-\mu\int_\Omega(\tau^n_-)^{2^{{m+1}}-1}\dive(u_t^n)\rd x.
\end{align}
We shall estimate the right-hand side of the above equality.  With the help of the Young inequality, the H\"older inequality, Sobolev embedding $H^1\hookrightarrow L^6$ and Lemma \ref{Lemappenergy}, we obtain
\begin{align}\label{tau-3}
	-\mu\int_\Omega&\dive(u_t^n)(\tau^n_-)^{2^{{m+1}}-1}\rd x\notag\\
	&\leqslant|\mu|\cdot\frac{2^{m+1}-1}{2^m}\int_\Omega u_t^n(\tau^n_-)^{2^m-1}\nabla\left((\tau^n_-)^{2^m}\right)\rd x\notag\\
	&\leqslant C\int_\Omega u_t^n\left(1^{2^m}+(\tau^n_-)^{(2^m-1)\cdot\frac{2^m}{2^m-1}}\right)\nabla\left((\tau^n_-)^{2^m}\right)\rd x\notag\\
	&\leqslant C\|u_t^n\|_{L^2(\Omega)}\|\nabla\left((\tau^n_-)^{2^m}\right)\|_{L^2(\Omega)}
+C\|u_t^n\|_{L^6(\Omega)}\|(\tau^n_-)^{2^m}\|_{L^3(\Omega)}\|\nabla\left((\tau^n_-)^{2^m}\right)\|_{L^2(\Omega)}\notag\\
	&\leqslant \frac{\alpha_m}{4}\|\nabla\left((\tau^n_-)^{2^m}\right)\|_{L^2(\Omega)}^2+C\alpha_m^{-1}
	+C\|\nabla u_t^n\|_{L^2(\Omega)}\|(\tau^n_-)^{2^m}\|_{L^3(\Omega)}\|\nabla\left((\tau^n_-)^{2^m}\right)\|_{L^2(\Omega)}\notag\\
	&\leqslant \frac{\alpha_m}{4}\|\nabla\left((\tau^n_-)^{2^m}\right)\|_{L^2(\Omega)}^2+C\alpha_m^{-1}
	+C\|(\tau^n_-)^{2^m}\|_{L^3(\Omega)}\|\nabla\left((\tau^n_-)^{2^m}\right)\|_{L^2(\Omega)},
\end{align}
where $\alpha_m:=\frac{2^{m+1}-1}{4^m}$.
The Gagliardo-Nirenberg inequality gives that
\begin{align}\label{tau-4}
	\|(\tau^n_-)^{2^m}\|_{L^3(\Omega)}\leqslant C\|\nabla\left((\tau^n_-)^{2^m}\right)\|_{L^2(\Omega)}^{\frac 45}\|(\tau^n_-)^{2^m}\|_{L^1(\Omega)}^{\frac 15}+C\|(\tau^n_-)^{2^m}\|_{L^1(\Omega)}.
\end{align}	
Substituting \eqref{tau-3} and \eqref{tau-4} into \eqref{tau-2} and using the Young inequality, we get	
\begin{align*}
	\frac{1}{2^{m+1}}\frac{\rd}{\rd t}\int_\Omega(\tau^n_-)^{2^{m+1}}\rd x
	&+\alpha_m\int_\Omega\left[\nabla\left((\tau^n_-)^{2^m}\right)\right]^2\rd x\\
	&\leqslant \frac{\alpha_m}{4}\|\nabla\left((\tau^n_-)^{2^m}\right)\|_{L^2(\Omega)}^2+C\alpha_m^{-1}
	+C\|\nabla\left((\tau^n_-)^{2^m}\right)\|_{L^2(\Omega)}^{\frac 95}\|(\tau^n_-)^{2^m}\|_{L^1(\Omega)}^{\frac 15}\\
	&~~~+C\|(\tau^n_-)^{2^m}\|_{L^1(\Omega)}\|\nabla\left((\tau^n_-)^{2^m}\right)\|_{L^2(\Omega)}\\
	&\leqslant \frac{\alpha_m}{2}\|\nabla\left((\tau^n_-)^{2^m}\right)\|_{L^2(\Omega)}^2+C\alpha_m^{-1}
	+C\left(\alpha_m^{-1}+\alpha_m^{-9}\right)\|(\tau^n_-)^{2^m}\|_{L^1(\Omega)}^{2},
\end{align*}	
which implies that the following inequality
\begin{align*}
 \frac{1}{2^{m+1}}\frac{\rd}{\rd t}\int_\Omega(\tau^n_-)^{2^{m+1}}\rd x
 +\frac{\alpha_m}{2}\int_\Omega\left[\nabla\left((\tau^n_-)^{2^m}\right)\right]^2\rd x
    \leqslant C\alpha_m^{-1}
    +C\left(\alpha_m^{-1}+\alpha_m^{-9}\right)\|(\tau^n_-)^{2^m}\|_{L^1(\Omega)}^{2},
\end{align*}	
is satisfied. Applying the Gagliardo-Nirenberg inequality again to the norm $\|\nabla\left((\tau^n_-)^{2^m}\right)\|_{L^2(\Omega)}^2$, we obtain	
\begin{align*}
	\frac{1}{2^{m+1}}\frac{\rd}{\rd t}\int_\Omega(\tau^n_-)^{2^{m+1}}\rd x
	+\frac{\alpha_m}{2}\int_\Omega(\tau^n_-)^{2^{m+1}}\rd x
	\leqslant C\alpha_m^{-1}+C\left(\alpha_m^{-1}+\alpha_m^{-9}+1\right)\|(\tau^n_-)^{2^m}\|_{L^1(\Omega)}^{2}.
\end{align*}
Then we can rewrite the above inequality in the following way
\begin{align*}
 \frac{\rd}{\rd t}\|\tau^n_-\|_{L^{2^{m+1}}(\Omega)}^{2^{m+1}}+\|\tau^n_-\|_{L^{2^{m+1}}(\Omega)}^{2^{m+1}}
  \leqslant C2^{10m}\left(\|\tau^n_-\|_{L^\infty\left(0,\infty; {L^{2^{m}}(\Omega)}\right)}^{2^{m+1}}+1\right).
\end{align*}	
It leads us to	
\begin{align*}
	\|\tau^n_-\|_{L^{2^{m+1}}(\Omega)}^{2^{m+1}}
	   &\leqslant e^{-t}\|\tau^n_-(0,\cdot)\|_{L^{2^{m+1}}(\Omega)}^{2^{m+1}}
	   +C2^{10m}\left(\|\tau^n_-\|_{L^\infty\left(0,\infty; {L^{2^{m}}(\Omega)}\right)}^{2^{m+1}}+1\right)\\
	   &\leqslant |\Omega|\cdot\|\tau^n_-(0,\cdot)\|_{L^{\infty}(\Omega)}^{2^{m+1}}
	   +C2^{10m}\left(\|\tau^n_-\|_{L^\infty\left(0,\infty; {L^{2^{m}}(\Omega)}\right)}^{2^{m+1}}+1\right),
\end{align*}
which indicates that
\begin{align*}
	\|\tau^n_-\|_{L^\infty\left(0,\infty; L^{2^{m+1}}(\Omega)\right)}
	\leqslant\left(|\Omega|\cdot\|\tau^n_-(0,\cdot)\|_{L^{\infty}(\Omega)}^{2^{m+1}}
	+C2^{10m}\left(\|\tau^n_-\|_{L^\infty\left(0,\infty; L^{2^{m}}(\Omega)\right)}^{2^{m+1}}+1\right)\right)^\frac{1}{2^{m+1}}.
\end{align*}

\textbf{Step 3. Start iteration.} Let us introduce a notation
$$\beta_m=\max\left\{\|\tau^n_-(0,\cdot)\|_{L^\infty(\Omega)},~\|\tau^n_-\|_{L^{\infty}\left(0,\infty; L^{2^m}(\Omega)\right)},~1\right\}.$$
Then, we arrive at
$$\beta_{m+1}\leqslant \beta_mC^{\frac{1}{2^{m+1}}}2^\frac{10m}{2^{m+1}},$$
which gives us that
$$\beta_m\leqslant \beta_1\prod_{k=2}^{m}C^{\frac{1}{2^{k}}}2^\frac{10(k-1)}{2^{k}}
=\beta_1C^{\sum_{k=2}^{m}{\frac{1}{2^{k}}}}2^{\sum_{k=2}^{m}\frac{10(k-1)}{2^{k}}}
\leqslant\beta_1C^{\sum_{k=2}^{\infty}{\frac{1}{2^{k}}}}2^{\sum_{k=2}^{\infty}\frac{10(k-1)}{2^{k}}}.$$
The convergence of both exponential series in the preceding inequality directly yields the uniform boundedness of the sequence $\beta_m$ for all $m\in \mathbb{N}$.
Therefore we obtain
$$\tau^n_-\in L^\infty\left((0,\infty)\times\Omega\right).$$

\textbf{Step 4. Take the limit of $\tau^n$.} By Proposition \ref{ProtauC}, we know that
$$\tau^n\in L^2\left(0,\infty; H^1(\Omega)\right),~~~\tau^n_t\in L^1\left(0,\infty; L^1(\Omega)\right).$$
Then, by the Aubin-Lions Lemma, we know that there exist a subsequence of $\tau^n$ and $\tau\in L^2\left(0,T; L^2(\Omega)\right)$ such that
$$\tau^n\rightarrow\tau~~~\text{in}~L^2\left(0,T; L^2(\Omega)\right).$$
It yields that (possibly on subsequence, see \cite{Brezis2011})
\begin{align}\label{taun1C}
	\tau^n(t,x)\rightarrow\tau(t,x),~~~\text{a.e.~in}~[0,T]\times\Omega.
\end{align}
On the other hand, \eqref{thetaL2C} gives us that (possibly on subsequence)
\begin{align}\label{theta1C}
	\theta^n(t,x)\rightarrow\theta(t,x),~~~\text{a.e.~in}~[0,T]\times\Omega.
\end{align}
Since $\tau^n=\log\theta^n$, thanks to \eqref{taun1C} and \eqref{theta1C},
we get $$\tau=\log\theta.$$
With the help of the Lemma \ref{Ththetainfty}, it leads us to
$$-\|\tau\|_{L^\infty\left((0,\infty)\times\Omega\right)}
\leqslant\tau\leqslant\|\tau\|_{L^\infty\left((0,\infty)\times\Omega\right)}
\leqslant\log\|\theta\|_{L^\infty\left((0,\infty)\times\Omega\right)}.$$
It gives us that the following inequality
$$0<\theta^*:=\re^{-\|\tau\|_{L^\infty\left((0,\infty)\times\Omega\right)}}
\leqslant\theta(t,x)\leqslant\|\theta\|_{L^\infty\left((0,\infty)\times\Omega\right)},$$
is satisfied for almost all $(t,x)\in[0,\infty)\times\Omega$.
Thus, the proof of Theorem \ref{Thtauinfty} is finished.	
\end{proof}

\subsection{Uniqueness of solutions}
In this section, we establish the uniqueness assertion of Theorem \ref{ThexistandPos} through the following theorem.
\begin{theorem}\label{ThUinque}
	The solution satisfying \eqref{equ1.1} in the sense of Definition \ref{solution} is unique.
\end{theorem}
\begin{proof}
Let us assume that $(u_1,\theta_1)$	and $(u_2,\theta_2)$ are two different solutions of \eqref{equ1.1} in the sense of Definition \ref{solution}.
We denote $u=u_1-u_2$ and $\theta=\theta_1-\theta_2$ and subtract the equations for $(u_2,\theta_2)$ from the equations for $(u_1,\theta_1)$.
It yields that
\begin{align*}
	\int_\Omega u_{tt}\varphi \rd x+\int_\Omega \nabla u\nabla\varphi \rd x+\int_\Omega\nabla u_{tt}\nabla\varphi\rd x
	=\mu \int_\Omega\nabla\theta\cdot\varphi\rd x.
\end{align*}
Taking $\varphi=u_t$ in the above equation, we get
\begin{align}\label{u12t}
\frac 12\frac{\rd}{\rd t}\left(\|u_t\|_{L^2(\Omega)}^2+\|\nabla u\|_{L^2(\Omega)}^2+\|\nabla u_t\|_{L^2(\Omega)}^2\right)
	\leqslant \frac14\|\nabla\theta\|_{L^2(\Omega)}^2+C\|u_t\|_{L^2(\Omega)}^2.
\end{align}
Similarly, we proceed with the heat equation and test it by $\psi=\theta$. Then we have
\begin{align}\label{I123}
	\frac 12\frac{\rd}{\rd t}\|\theta\|_{L^2(\Omega)}^2+\|\nabla\theta\|_{L^2(\Omega)}^2
	&=\mu\int_\Omega\theta^2\dive(u_{1,t})\rd x+\mu\int_\Omega\theta_2\theta\dive(u_t)\rd x\notag\\
	&\leqslant-2\mu\int_\Omega u_{1,t}\theta\nabla\theta\rd x
	-\mu\int_\Omega\theta_2\nabla\theta u_t\rd x-\mu\int_\Omega\theta\nabla\theta_2u_t\rd x\notag\\
	&:=I_1+I_2+I_3.
\end{align}	
Now we shall estimate the terms $I_1,~I_2$ and $I_3$. With the help of the H\"older inequality, Sobolev embedding $H^1\hookrightarrow L^6$, \eqref{energy} and the Young inequality, we obtain
\begin{align*}
	I_1&=-2\mu\int_\Omega u_{1,t}\theta\nabla\theta\rd x
	\leqslant2|\mu|\cdot\|u_{1,t}\|_{L^6(\Omega)}\|\nabla\theta\|_{L^2(\Omega)}\|\theta\|_{L^3(\Omega)}\\
	&\leqslant C\|\nabla u_{1,t}\|_{L^2(\Omega)}
	\|\nabla\theta\|_{L^2(\Omega)}^{\frac 32}\|\theta\|_{L^2(\Omega)}^{\frac 12}+C\|\nabla u_{1,t}\|_{L^2(\Omega)}	\|\nabla\theta\|_{L^2(\Omega)}\|\theta\|_{L^2(\Omega)}\\
	&\leqslant C
	\|\nabla\theta\|_{L^2(\Omega)}^{\frac 32}\|\theta\|_{L^2(\Omega)}^{\frac 12}+C\|\nabla\theta\|_{L^2(\Omega)}\|\theta\|_{L^2(\Omega)}
	\leqslant\frac14\|\nabla\theta\|_{L^2(\Omega)}^2+C\|\theta\|_{L^2(\Omega)}^2,
\end{align*}	
Here we have used the following Gagliardo-Nirenberg inequality in the above second inequality 
\begin{align}\label{GN322}
	\|\theta\|_{L^3(\Omega)}\leqslant C\|\nabla\theta\|_{L^2(\Omega)}^{\frac 12}\|\theta\|_{L^2(\Omega)}^{\frac 12}+C\|\theta\|_{L^2(\Omega)}.
\end{align}
With $I_2$, using the H\"older inequality, Lemma \ref{Ththetainfty} and the Young inequality, we have
\begin{align*}
	I_2&=-\mu\int_\Omega\theta_2\nabla\theta u_t\rd x\\
	&\leqslant |\mu|\cdot\|\theta_2\|_{L^\infty\left((0,\infty)\times\Omega\right)}
	\|\nabla\theta\|_{L^2(\Omega)}\|u_t\|_{L^2(\Omega)}\\
	&\leqslant\frac14\|\nabla\theta\|_{L^2(\Omega)}^2+C\|u_t\|_{L^2(\Omega)}^2.
\end{align*}	
Next, we will devote to estimate $I_3$. Applying the H\"older inequality, Sobolev embedding $H^1\hookrightarrow L^6$, Gagliardo-Nirenberg inequality \eqref{GN322}, the Young inequality again and Theorem \ref{Thhigher}, we get
\begin{align*}
	I_3&=-\mu\int_\Omega\theta\nabla\theta_2u_t\rd x\\
	&\leqslant|\mu|\cdot\|\theta\|_{L^3(\Omega)}\|\nabla\theta_2\|_{L^2(\Omega)}\|u_t\|_{L^6(\Omega)}
	\leqslant C\|\theta\|_{L^3(\Omega)}\|\nabla u_t\|_{L^2(\Omega)}\\
	&\leqslant C\|\nabla\theta\|_{L^2(\Omega)}^{\frac 12}\|\theta\|_{L^2(\Omega)}^{\frac 12}\|\nabla u_t\|_{L^2(\Omega)}+C\|\theta\|_{L^2(\Omega)}\|\nabla u_t\|_{L^2(\Omega)}\\
	&\leqslant C\|\nabla\theta\|_{L^2(\Omega)}\|\theta\|_{L^2(\Omega)}+C\|\nabla u_t\|_{L^2(\Omega)}^2+C\|\theta\|_{L^2(\Omega)}\|\nabla u_t\|_{L^2(\Omega)}\\
	 &\leqslant \frac14\|\nabla\theta\|_{L^2(\Omega)}^2+C\|\theta\|_{L^2(\Omega)}^2
	 +C\|\nabla u_t\|_{L^2(\Omega)}^2.
\end{align*}
Plugging the above estimates of $I_1,~I_2$ and $I_3$ into \eqref{I123}, it yields that
\begin{align*}
	\frac 12\frac{\rd}{\rd t}\|\theta\|_{L^2(\Omega)}^2+\frac14\|\nabla\theta\|_{L^2(\Omega)}^2
	\leqslant C\left(\|\theta\|_{L^2(\Omega)}^2+\|u_t\|_{L^2(\Omega)}^2+\|\nabla u_t\|_{L^2(\Omega)}^2\right).
\end{align*}
Adding the above inequality to \eqref{u12t}, we arrive at
\begin{align*}
	\frac{\rd}{\rd t}\bigg(\|u_t\|_{L^2(\Omega)}^2&+\|\nabla u\|_{L^2(\Omega)}^2
	+\|\nabla u_t\|_{L^2(\Omega)}^2+\|\theta\|_{L^2(\Omega)}^2\bigg)\\
	&\leqslant C\left(\|u_t\|_{L^2(\Omega)}^2+\|\nabla u\|_{L^2(\Omega)}^2
	+\|\nabla u_t\|_{L^2(\Omega)}^2+\|\theta\|_{L^2(\Omega)}^2\right).
\end{align*}
Thus, by using the Gronwall's inequality, we obtain
\begin{align*}
	\|u_t\|_{L^2(\Omega)}^2&+\|\nabla u\|_{L^2(\Omega)}^2
	+\|\nabla u_t\|_{L^2(\Omega)}^2+\|\theta\|_{L^2(\Omega)}^2
	\leqslant0,
\end{align*}
which implies $u_t=\nabla u=\nabla u_t=\theta=0$, for almost all $t\in[0,T]$ and any $T>0$.
Therfore we get $u_1=u_2$ and $\theta_1=\theta_2$. 
\end{proof}
\section{Asymptotic behavior}
In this section, we establish the proof of Theorem \ref{Th-Convergences}, which characterizes the asymptotic behavior of the heat-conducting elastic body.
We first establish the stability for \eqref{equ1.1} with respect to initial data in proper functional spaces, as detailed in Proposition \ref{Prostability}. Subsequently, we construct the dynamical system and its associated phase space describing the solution evolution of \eqref{equ1.1}. This framework enables the identification of $\omega$-limit sets.
Proposition \ref{ProtauC} provides a quantitative formulation of the second law of thermodynamics, which serves as a cornerstone of our analysis.
More precisely, our proof shows that the asymptotic states of the heat-conducting elastic body are identified by the second law of thermodynamics.

Now, we present the stability of \eqref{equ1.1} in proper functional spaces with respect to initial data.
\begin{proposition}\label{Prostability}
Let $(u_{0,1}, v_{0,1}, \theta_{0,1})$ and $(u_{0,2}, v_{0,2}, \theta_{0,2})$ be two triples of initial data such that
$$ u_{0,i}, v_{0,i} \in H^2(\Omega) \cap H_0^1(\Omega), \quad \theta_{0,i} \in H^1(\Omega), \quad i=1,2. $$
Let us also assume that $(u_1,\theta_1)$ and $(u_2,\theta_2)$ are the weak solution of \eqref{equ1.1} starting from $(u_{0,1},v_{0,1},\theta_{0,1})$ and $(u_{0,2},v_{0,2},\theta_{0,2})$, respectively.
Then, we obtain that
\begin{align*}
	\|u_1-u_2\|_{H^1_0(\Omega)}&+\|u_{1,t}-u_{2,t}\|_{H_0^1(\Omega)}+\|\theta_1-\theta_2\|_{L^2(\Omega)}\\
	&\leqslant C\left(\|u_{0,1}-u_{0,2}\|_{H^1_0(\Omega)}+\|v_{0,1}-v_{0,2}\|_{L^2(\Omega)}+\|\theta_{0,1}-\theta_{0,2}\|_{L^2(\Omega)}\right).
\end{align*}
\end{proposition}
\begin{proof}
Let us denote $u=u_1-u_2$ and $\theta=\theta_1-\theta_2$ and subtract the equations for $(u_2,\theta_2)$ from the equations for $(u_1,\theta_1)$.
In a manner analogous to the proof steps of Theorem \ref{ThUinque}, we can derive
\begin{align*}
	\frac{\rd}{\rd t}\bigg(\|u_t\|_{L^2(\Omega)}^2&+\|\nabla u\|_{L^2(\Omega)}^2
	+\|\nabla u_t\|_{L^2(\Omega)}^2+\|\theta\|_{L^2(\Omega)}^2\bigg)\\
	&\leqslant C\left(\|u_t\|_{L^2(\Omega)}^2+\|\nabla u\|_{L^2(\Omega)}^2
	+\|\nabla u_t\|_{L^2(\Omega)}^2+\|\theta\|_{L^2(\Omega)}^2\right).
\end{align*}
Applying the Gronwall's inequality, we complete the proof of Proposition \ref{Prostability}.	
\end{proof}

Next, we construct the corresponding phase space and dynamical system, which will be crucial for proving Theorem \ref{Th-Convergences}.
We introduce the following notation
$$\mathcal{M}:=\left\{\Theta\in H^1(\Omega);~{\rm ess}\inf\limits_{x\in\Omega}\Theta(x)>0\right\}.$$
We define an operator
\begin{align*}
	S(t):&\left(H_0^1(\Omega)\cap H^2(\Omega)\right)\times\left(H_0^1(\Omega)\cap H^2(\Omega)\right)\times \mathcal{M}\\
	&\rightarrow\left(H_0^1(\Omega)\cap H^2(\Omega)\right)\times\left(H_0^1(\Omega)\cap H^2(\Omega)\right)\times \mathcal{M},
\end{align*}
by formula
$$S(t)\left(u_0,v_0,\theta_0\right)=\left(u_t(t),v_t(t),\theta_t(t)\right),$$
where $u$ and $\theta$ be a solution of \eqref{equ1.1} starting from initial data $u_0,~v_0$ and $\theta_0$.
With the help of Theorem \ref{ThUinque}, for all $t, h\in[0,\infty)$, we arrive at
$$S(t+h)=S(t)\circ S(h).$$

Proposition \ref{Prostability} gives us the following significant corollary.
\begin{corollary}\label{ColoSt}
Let
$(u_0,v_0,\theta_0),~(u^n,v^n,\theta^n)\in
		\left(H_0^1(\Omega)\cap H^2(\Omega)\right)\times\left(H_0^1(\Omega)\cap H^2(\Omega)\right)\times \mathcal{M}$
such that
$$(u_0,v_0,\theta_0)\rightarrow(u^n,v^n,\theta^n)~~\text{in}~H_0^1(\Omega)\times H_0^1(\Omega)\times L^2(\Omega).$$
Then, for each $t\geqslant0$, we also have
$$S((t))(u_0,v_0,\theta_0)\rightarrow S(t)(u^n,v^n,\theta^n)~~\text{in}~H_0^1(\Omega)\times H_0^1(\Omega)\times L^2(\Omega).$$
\end{corollary}

Now, we go back to the proof of Theorem \ref{Th-Convergences}.

\textbf{Proof of Theorem \ref{Th-Convergences}.}
Let us take a measurable set $X\subset [0,\infty)$ such that $meas([0,\infty)\setminus X) = 0$.
By Theorem \ref{Thexist1}, we have a subsequence $\{t_k\}\subset X$ (still denoted as $\{t_k\}$) satisfies
\begin{align*}
	u(t_k,\cdot)&\rightharpoonup\tilde{u}~~\text{in}~H_0^1(\Omega)\cap H^2(\Omega),\\
	u_t(t_k,\cdot)&\rightharpoonup\tilde{v}~~\text{in}~H_0^1(\Omega)\cap H^2(\Omega),\\
	\theta(t_k,\cdot)&\rightharpoonup\tilde{\theta}~~\text{in}~H^1(\Omega),
\end{align*}
where $\tilde{u}\in H_0^1(\Omega)\cap H^2(\Omega)$,
$\tilde{v}\in H_0^1(\Omega)\cap H^2(\Omega)$ and $\tilde\theta\in H^1(\Omega)$ are certain functions.
By the compact embedding $H^2\hookrightarrow\hookrightarrow H_0^1,~H^1\hookrightarrow\hookrightarrow L^2$ and \eqref{C0T}, we get
\begin{align*}
	u(t_k,\cdot)&\rightarrow\tilde{u}~~\text{in}~H_0^1(\Omega),\\
	u_t(t_k,\cdot)&\rightarrow\tilde{v}~~\text{in}~H_0^1(\Omega),\\
	\theta(t_k,\cdot)&\rightarrow\tilde{\theta}~~\text{in}~L^2(\Omega).
\end{align*}
Let us take $h>0$ and consider the sequences $u(t_k+h,\cdot),~u_t(t_k+h,\cdot)$ and $\theta(t_k+h,\cdot)$.
In a similar manner to the above case, there exist $\hat{u}\in H_0^1(\Omega)\cap H^2(\Omega)$,
$\hat{v}\in H_0^1(\Omega)\cap H^2(\Omega)$ and $\hat\theta\in H^1(\Omega)$
such that (possibly on subsequence of $\{t_k\}$)
\begin{align*}
	u(t_k+h,\cdot)&\rightarrow\hat{u}~~\text{in}~H_0^1(\Omega),\\
	u_t(t_k+h,\cdot)&\rightarrow\hat{v}~~\text{in}~H_0^1(\Omega),\\
	\theta(t_k+h,\cdot)&\rightarrow\hat{\theta}~~\text{in}~L^2(\Omega).
\end{align*}
With the help of Corollary \ref{ColoSt}, we obtain that
\begin{align}\label{u-bar}
	(\hat u, \hat v, \hat\theta)
	&=\lim\limits_{k\rightarrow\infty}\left(u(t_k+h,\cdot),u_t(t_k+h,\cdot),\theta(t_k+h,\cdot)\right)\notag\\
	&=\lim\limits_{k\rightarrow\infty}S(h)\left(u(t_k,\cdot),u_t(t_k,\cdot),\theta(t_k,\cdot)\right)\notag\\
	&=S(h)\left(\tilde u,\tilde{v},\tilde\theta\right)
	:=\left(\bar{u}(h,\cdot),\bar{v}(h,\cdot),\bar\theta(h,\cdot)\right).
\end{align}
Here we denoted by $\bar{u}$ and $\bar\theta$ the weak solution of \eqref{equ1.1} with the initial values $\tilde u,\tilde{v}$ and $\tilde\theta$ (they can be established because the limiting objects possess sufficient regularity to allow for the definition of weak solutions).

On the other hand, Proposition \ref{ProtauC} ensures that the following function
$$t\mapsto\int_\Omega\log\theta(t,\cdot)\rd x,$$
is non-decreasing and bounded. According to Theorem \ref{Thtauinfty}, the Lebesgue Dominated Convergence Theorem and the uniqueness of the limit, we get
\begin{align*}
	\int_\Omega\log\tilde\theta\rd x
	&=\lim\limits_{k\rightarrow\infty}\int_\Omega\log\theta(t_k,\cdot)\rd x\\
	&=\lim\limits_{k\rightarrow\infty}\int_\Omega\log\theta(t_k+h,\cdot)\rd x
	=\int_\Omega\log\bar\theta(h,\cdot)\rd x,
\end{align*}
which implies that $\int_\Omega\log\bar{\theta}(t,\cdot)\rd x$ is constant at the $\omega$-limit set.

Now, we claim that $\bar\theta$ is constant in space and time.
By Proposition \ref{ProtauC}, we obtain
\begin{align*}
	0=\frac{\rd}{\rd t}\int_\Omega\log\bar{\theta}\rd x=\int_\Omega|\nabla(\log\bar{\theta})|^2\rd x.
\end{align*}
Thus, we have $\nabla(\log\bar{\theta})=0$, for any $t>0$ and $x\in\Omega$. So, for each fixed $t>0$, $\bar{\theta}$ is constant in space. Further, applying Proposition \ref{ProtauC} again, we obtain
\begin{align*}
	0=\frac{\rd}{\rd t}\int_\Omega\log\bar{\theta}\rd x
	=|\Omega|\cdot\frac{\rd}{\rd t}\left(\log\bar{\theta}\right),
\end{align*}
which means $\bar{\theta}$ is constant in time.
Thus, $\bar{\theta}$ is constant in space and time.

So far, we have identified all the trajectories of \eqref{equ1.1} that originate from the functions $\tilde{\theta}$. These functions $\tilde{\theta}$ represent all the potential limits over time of the temperature satisfying \eqref{equ1.1}.
Assume that there exist two different functions $\tilde{\theta}_1$ and $\tilde{\theta}_2$,
which are limits of $\theta(t_{km})$ and $\theta(t_{kl})$ for different time sequences $t_{km},~t_{kl}\rightarrow\infty$, respectively. By utilizing Proposition \ref{ProtauC} again, we arrive at
\begin{align*}
	\int_\Omega\log\tilde{\theta}_1\rd x=\int_\Omega\log\tilde{\theta}_2\rd x.
\end{align*}
Then we get $\tilde{\theta}_1=\tilde{\theta}_2$, since we already know that $\tilde{\theta}_1$ and $\tilde{\theta}_2$ are constant in space and time. Thus, we notice that the $\omega$-limit set of a solution to \eqref{equ1.1} actually includes only one potential limit for temperature $\tilde{\theta}$.

Since $\bar{\theta}$ is constant and positive (by Theorem \ref{Thtauinfty}), we get
\begin{equation}\label{bar0}
	\begin{cases}
		\bar u_{tt}-\Delta\bar u-\Delta\bar u_{tt}=0,~&(x,t)\in\Omega\times(0,T),\\
		0=\mu \bar\theta\dive(u_t),~&(x,t)\in\Omega\times(0,T),\\
		\bar u=0,~&(x,t)\in\partial\Omega\times(0,T),\\
	\end{cases}
\end{equation}
where $\bar{u}$ and $\bar\theta$ the weak solution of \eqref{equ1.1} with the initial values $\tilde u,\tilde{v}$ and $\tilde\theta$. Combining the second equation and boundary conditions in \eqref{bar0}, we get $\bar u_t=0$ for all $(t,x)\in[0,\infty)\times\Omega$, which means that $\bar{u}$ is constant in time for each fixed $x\in\Omega$. By utilizing the first equation in \eqref{bar0} and boundary conditions again, we obtain $\bar{u}=0$. Together with \eqref{u-bar}, we get
\begin{align*}
	u(t_k+h,\cdot)&\rightarrow0~~\text{in}~H_0^1(\Omega),\\
	u_t(t_k+h,\cdot)&\rightarrow0~~\text{in}~H_0^1(\Omega).
\end{align*}
Therefore, for a measurable set $X \subset [0,\infty)$ with $meas([0,\infty)\setminus X) = 0$, the limit as $t \to \infty$ (where $t \in X$) yields that
\begin{align*}
	\lim_{\substack{t\rightarrow\infty \\ t \in X}} u(t,\cdot)= 0~\text{in}~H_0^1(\Omega)
	~~~\text{and}~~~\lim_{\substack{t\rightarrow\infty \\ t \in X}} u(t,\cdot)= 0~\text{in}~H_0^1(\Omega).
\end{align*}
Then, by Theorem \ref{ThexistandPos}, we obtain
\begin{align*}
	\lim_{\substack{t\rightarrow\infty}} u(t,\cdot)= 0~\text{in}~H_0^1(\Omega)
	~~~\text{and}~~~\lim_{\substack{t\rightarrow\infty}} u(t,\cdot)= 0~\text{in}~H_0^1(\Omega).
\end{align*}
Finally, together with the above the limits and Proposition \ref{energy}, we deduce that
 \begin{align*}
 	&~~~\frac 12\int_\Omega v_0^2\rd x+\frac 12\int_\Omega|\nabla u_{0}|^2\rd x+\frac{1}{2}\int_\Omega|\nabla v_0|^2\rd x+\int_\Omega\theta_0 \rd x\\
 	&=\lim\limits_{k\rightarrow\infty}\left(\frac 12\int_\Omega u_t^2(t_k,\cdot)\rd x+\frac 12\int_\Omega |\nabla u(t_k,\cdot)|^2\rd x+\frac{1}{2}\int_\Omega|\nabla u_t(t_k,\cdot)|^2\rd x+\int_\Omega\theta(t_k,\cdot) \rd x\right)\\
 	&=\int_\Omega\tilde{\theta}\rd x=\tilde{\theta}\cdot|\Omega|.
 \end{align*}
Therefore, we get
$$\theta_\infty:=\tilde\theta=\left(\frac 12\int_\Omega v_0^2\rd x+\frac 12\int_\Omega|\nabla u_{0}|^2\rd x+\frac{1}{2}\int_\Omega|\nabla v_0|^2\rd x+\int_\Omega\theta_0 \rd x\right)\cdot|\Omega|^{-1},$$
which implies
$$\lim_{\substack{t\rightarrow\infty}} \theta(t,\cdot)=\theta_\infty~\text{in}~L^2(\Omega).$$
The proof of Theorem \ref{Th-Convergences} is now completed.
$\hfill\square$

	
\end{document}